\documentclass[12pt,a4paper]{article}
\usepackage[cp1251]{inputenc}
\usepackage{hyperref}
\usepackage[english]{babel}
\usepackage{amsmath,amsthm,dsfont,pifont,amsfonts,MnSymbol}
\usepackage{MPatr_Style_English_02.2020}

\author{Mikhail Patrakeev\footnote{Krasovskii Institute of Mathematics and Mechanics of UB RAS, 620108, 16 Sofia Kovalevskaya street, Yekaterinburg, Russia; \textit{e-mail address}\textup{:} patrakeev@mail.ru}
\footnote{This work is a part of research conducted in the Ural Mathematical Center}}

\title{Reverse induction proof of D property \\ of the countable power of the Sorgenfrey line}

\date{}

\begin{document}
\mathsurround=1pt
\hyphenation{tran-si-ti-ve par-ti-al Sor-gen-frey neigh-bour-hood }
\renewcommand{\proofname}{\textup{\textbf{Proof}}}
\renewcommand{\abstractname}{\textup{Abstract}}
\renewcommand{\refname}{\textup{References}}
\newcommand{\nos}{\mathsurround=0pt}
\renewcommand{\iff}{\quad{\colon}{\longleftrightarrow}\quad}
\renewcommand{\SS}{\mathbb{S}}
\newcommand{\ZZ}{\mathbb{Z}}
\newcommand{\broom}{\mathsf{broom}}
\newcommand{\cone}{\mathsf{cone}}
\newcommand{\crown}{\mathsf{crown}}
\newcommand{\edge}{\mathsf{edge}}
\newcommand{\bottom}{\mathsf{bottom}}
\newcommand{\plane}{\mathsf{plane}}
\newcommand{\facet}{\mathsf{facet}}
\newcommand{\hcone}{\breve\cone}
\newcommand{\hcrown}{\breve\crown}
\newcommand{\hedge}{\breve\edge}
\newcommand{\hbottom}{\breve\bottom}
\newcommand{\hplane}{\breve\plane}
\newcommand{\hfacet}{\breve\facet}
\renewcommand{\ll}{\langle}
\newcommand{\rr}{\rangle}
\newcommand{\res}{\,{\upharpoonright}\,}

\maketitle
\overfullrule5pt

\begin{abstract}
We present a new method of proof, which we call reverse induction. This method allows to establish certain properties of a product $\prod_{{i}=0}^{\infty}{X}_{i}$ by making a kind of ``reverse induction step'' from $\prod_{{i}={n}+1}^{\infty}{X}_{i}$ to
$\prod_{{i}={n}}^{\infty}{X}_{i}$ for an arbitrary natural ${n}$.
Using this method we answer a question posed by E.\,K.\,van Douwen and W.\,F.\,Pfeffer in 1979~\cite{pacific} by proving that the countable power of the Sorgenfrey line is a D-space.
\end{abstract}

\section{Introduction}


A topological space ${X}$ is a \emph{D-space} iff whenever $\varphi({p})$ is a open neighbourhood of~${p}$ for each ${p}\in{X}$, there is a closed discrete\footnote{A set ${D}$ is \emph{discrete} iff each point ${q}\in{D}$ has an open neighbourhood ${O}({q})$ such that ${O}({q})\cap{D}=\{{q}\}$} subset ${D}$ of ${X}$ such that $\bigcup_{{p}\in{D}}\varphi({p})={X}$.
This notion was introduced by E.\,K.\,van Douwen and W.\,F.\,Pfeffer in~\cite{pacific}. Every ${T}_1\mathsurround=0pt$ compact space is a D-space, so they asked weather the D property follows from covering properties (see~\cite{handbook}) such as the Lindel\"{o}ff property, paracompactness, subparacompactness, and so on. Most of these questions are still open although nowadays they are central in set-theoretic topology, see~\cite{D_open.pr,D}. 

The countable power of the Sorgenfrey line\footnote{The Sorgenfrey line is the real line with the topology generated by the base $\{[{a},{b}):{a}<{b}\}$} is a hereditarily subparacompact space~\cite{S.hered.subpar}. Pfeffer and van Douwen expressed a hope that this space could be a source of counterexamples to the above questions, so they asked~\cite{pacific} whether (subsets of) at most countable powers of the Sorgenfrey line possess the D property. They proved~\cite{pacific} that all finite powers of the Sorgenfrey line are D-spaces, and two years later, in 1981, Peter de Coux proved~\cite{deCaux} that all subsets of finite powers of the Sorgenfrey line also possess the D property.

In this paper we prove that the countable power of the Sorgenfrey line is also a D-space, see Theorem~\ref{teo.S.D}.
We receive this result by seemingly a new method of proof, which we call the principle of reverse induction, see Section~\ref{sect.RI}. This method allows to establish different properties similar to compactness for countable products of spaces. We formulate and prove this principle in the case of D property, see Theorem~\ref{teo.rev.ind}.
In Remark~\ref{rem.rev.ind} we show how to generalize this principle in a way to establish other properties such as, for example, different covering properties. We hope that this method will be used in other areas of mathematics to establish properties similar to compactness for products of other structures.

\section{Reverse induction}\label{sect.RI}

We use terminology from~\cite{top.enc} and~\cite{kun}. A \textit{space} is a topological space, $\omega$ is the set of natural numbers, $\mathsf{closure}({S},{X})$ is the closure of a set ${S}$ in a space ${X}$.
By ${f}\res{S}$ we denote the restriction a function ${f}$ to a set ${S}$.
So, if ${q}=\ll{q}_{i}\rr_{{i}\in{A}}\in\prod_{{i}\in{A}}{X}_{i}$
and ${B}\subseteq{A}$,
then ${q}\res{B}=\ll{q}_{i}\rr_{{i}\in{B}}\in\prod_{{i}\in{B}}{X}_{i}$;
in particular, ${q}\res\varnothing$ is the unique function with the empty domain.

\begin{defi}
Suppose that ${X}=\prod_{{i}\in\omega}{X}_{i}$ and ${A}\subseteq\omega$. Then

  \begin{itemize}
  \item [\ding{46}\:]
    ${P}\nos\,$ is an $\nos{A}$\emph{-broom} on ${X}\iff{P}\subseteq{X}$ and
    ${p}\res{A}={q}\res{A}$ for all ${p},{q}\in{P}$.
  \end{itemize}
\end{defi}

\noindent Note that a nonempty set ${P}$ is an $\nos{A}$-broom on $\prod_{{i}\in\omega}{X}_{i}$ \ iff \,\ there is ${x}\in\prod_{{i}\in{A}}{X}_{i}$ such that ${p}\res{A}={x}$ for all ${p}\in{P}$.

\begin{rema}\label{rem.brooms}
Suppose that ${X}=\prod_{{i}\in\omega}{X}_{i}$ is a product of $T_1\nos$-spaces,
${A},{B}\subseteq\omega$, and ${P}$ is an $\nos{A}$-broom on ${X}$.
Then:

  \begin{itemize}
  \item [\textup{a.}\,]
    $\mathsf{closure}({P},{X})\nos\,$ is an $\nos{A}$-broom on ${X}$.
  \item [\textup{b.}\,]
    If ${B}\subseteq{A}$,
    then ${P}$ is a $\nos{B}$-broom on ${X}$.
  \item [\textup{c.}\,]
    If ${R}\subseteq{P}$,
    then ${R}$ is an $\nos{A}$-broom on ${X}$.
  \item [\textup{d.}\,]
    If ${Q}$ is an $\nos{A}$-broom on ${X}$ and is a $\nos{B}$-broom on ${X}$,
    then ${Q}$ is an $\nos({A}\cup{B})$-broom on ${X}$.\hfill$\qed$%
  \end{itemize}
\end{rema}

Recall that $\varphi$ is an \emph{NA} (\emph{neighbourhood assignement})
for a space ${X}$ \,iff\,
$\varphi$ is a function from ${X}$ to the powerset of ${X}$
such that $\varphi({p})$ is an open neighbourhood of ${p}$ for all ${p}\in{X}$.
A space ${X}$ is called a {\it D-space} \,iff\, for every NA $\varphi$ for ${X}$
there exists a closed discrete set ${D}$ in ${X}$ such that
$\bigcup\{\varphi({p}):{p}\in{D}\}\supseteq{X}$.

\begin{nota}
Suppose that $\varphi$ is an NA for a space ${X}$ and ${Y}\subseteq{X}$. Then:

  \begin{itemize}
  \item [\ding{46}\:]
    $\varphi\ll{Y}\rr\:\coloneq\:\bigcup\{\varphi({p}):{p}\in{Y}\}\nos\,$.
  \item [\ding{46}\:]
    $\varphi\nos\,$ is \emph{correct} on ${Y}\iff$

    there exists a closed discrete set ${D}$ in the subspace ${Y}$ such that $\varphi\ll{D}\rr\supseteq{Y}$.
  \end{itemize}
\end{nota}

\noindent
Note that a space ${X}$ is a D-space \,iff\, every NA $\varphi$ for ${X}$ is correct on ${X}$.
\smallskip

In the following notation we assume that the natural numbers are identified with the finite ordinals, so each natural ${n}$ equals the set $\{0,1,\ldots,{n}-1\}$ of its predecessors; in particular, $0=\varnothing$. ``${P}\nos\,$ is a closed $\nos{n}$-broom on ${X}$'' means that ${P}$ is a closed subset of ${X}$ and ${P}$ is an $\nos{n}$-broom on ${X}$.
RIH stands for Reverse Induction Hypothesis.

\begin{nota}\label{nota.RIH}
Suppose that ${X}=\prod_{{i}\in\omega}{X}_{i}$ is a product of spaces, $\varphi$ is an NA for ${X}$, ${n}\in\omega$, and ${P}$ is a closed $\nos{n}$-broom on ${X}$. Then

  \begin{itemize}
  \item [\ding{46}\:]
    $\mathsf{RIH}({X},\varphi,{n},{P})\iff\nos\,$

    for every natural ${m}>{n}$ and every closed $\nos{m}$-broom ${Q}\subseteq{P}$ on ${X}$ $[\:\varphi$ is  correct on ${Q}\:]$.
  \end{itemize}
\end{nota}

\begin{teor}[\textbf{Reverse induction for D property in countable products}]\label{teo.rev.ind}\mbox{ }

\noindent Let ${X}=\prod_{{i}\in\omega}{X}_{i}$ be a product of ${T}_1\nos$-spaces.
Suppose that for every NA $\varphi$ for ${X}$, every ${n}\in\omega$, and every closed $\nos{n}$-broom ${P}$ on ${X}$,

  \begin{itemize}
  \item [\ding{70}\,]
    $\mathsf{RIH}({X},\varphi,{n},{P})\nos\,$
    implies $[\:\varphi$ is  correct on ${P}\:]$.
  \end{itemize}
Then ${X}$ is a D-space.
\end{teor}

\begin{proof}
Suppose on the contrary that ${X}$ is not a D-space.
Then there is an NA $\varphi$ for ${X}$ that is not correct on ${X}$.
Put ${n}_{0}\coloneq{0}$ and ${P}_{0}\coloneq{X}$. Then
${P}_{0}$ is a closed $\nos{n}_{0}$-broom on ${X}$ and $\varphi$ is not  correct on ${P}_{0}$.
It follows from (\ding{70}) that
$\mathsf{RIH}({X},\varphi,{n}_{0},{P}_{0})$ does not hold,
so there are natural ${n}_{1}>{n}_{0}$
and a closed $\nos{n}_{1}$-broom ${P}_{1}\subseteq{P}_{0}$ on ${X}$
such that $\varphi$ is not  correct on ${P}_{1}$.
Repeating this argument, it is straightforward to build, by recursion on ${j}$,
a sequence $\ll{n}_{j}\rr_{{j}\in\omega}$ of natural numbers and a sequence $\ll{P}_{j}\rr_{{j}\in\omega}$ of subsets of ${X}$ such that, for all ${j}\in\omega$,

  \begin{itemize}
  \item [1.\,]
    ${P}_{j}\nos\,$ is a closed $\nos{n}_{j}$-broom on ${X}\text{,}$
  \item [2.\,]
    $\varphi\nos\,$ is not  correct on ${P}_{j}$,
  \item [3.\,]
    ${n}_{{j}+1}>{n}_{j}\nos\,$,\quad and
  \item [4.\,]
    ${P}_{{j}+1}\subseteq{P}_{j}\nos\,$.
  \end{itemize}

It follows from (2) that each ${P}_{j}$ is not empty.
Then (1) implies that, for each ${j}\in\omega$,
there is ${x}_{j}\in\prod_{{i}\in{n}_{j}}{X}_{i}$
such that ${p}\res{{n}_{j}}={x}_{j}$ for all ${p}\in{P}_{j}$.
Also it follows from (4) that
${x}_{{j}+1}\res{{n}_{j}}={x}_{j}$ for all ${j}\in\omega$.
Then (3) implies that there is ${r}=\ll{r}_{i}\rr_{{i}\in\omega}\in{X}$
such that ${r}\res{{n}_{j}}={x}_{j}$ for all ${j}\in\omega$.
Let ${U}=\prod_{{i}\in{k}}{U}_{i}\times\prod_{{i}\in\omega\setminus{k}}{X}_{i}$
be a basic neighbourhood of ${r}$ in ${X}$
such that ${r}\in{U}\subseteq\varphi({r})$.
It follows from (3) that
there is ${j}'\in\omega$
such that ${n}_{{j}'}\geqslant{k}$,
so $\varphi({r})\supseteq{P}_{{j}'}$.
Using (3) and (4) it is not hard to show that ${r}\in\mathsf{closure}({P}_{{j}'},{X})={P}_{{j}'}$,
so $\{{r}\}$ is a closed discrete set in the subspace ${P}_{{j}'}$ because ${X}$ is a $T_1\nos$-space (the $T_1$ restriction is needed only here).
Therefore $\varphi\ll\{{r}\}\rr\supseteq{P}_{{j}'}$ contradicts (2).
\end{proof}
\smallskip

\begin{rema}\label{rem.rev.ind}\mbox{ }

The complete analogues of Theorem~\ref{teo.rev.ind} can be easily formulated and proved for many other properties similar to compactness, such as the Lindel\"{o}ff property and many other covering properties~\cite{top.enc,handbook}.

For example, if we want to get the compactness property variant, we are to define an open cover $\varphi$ of a space ${X}$ to be \emph{correct} on ${Y}\subseteq{X}$ \,iff\, there is a finite $\psi\subseteq\varphi$ such that $\bigcup\psi\supseteq{Y}$, replace the words ``NA for ${X}$'' by the words ``open cover of ${X}$'' in Notation~\ref{nota.RIH} and Theorem~\ref{teo.rev.ind}, remove the ${T}_1$ restriction, and not to forget to replace ``D'' with ``compact''. After that we obtain \emph{the principle of Reverse induction for compactness in countable products.}

The proof of such a principle remains the same. All we need from the compactness property  for this proof to work is that an open cover $\varphi$ of a space ${X}$ is correct on ${Y}\subseteq{X}$ whenever ${Y}=\varnothing$ or ${Y}\subseteq{U}$ for some ${U}\in\varphi$.
\end{rema}

As a trivial example of application of \emph{the principle of Reverse induction for compactness in countable products,} it is a simple exercise to show, by using it, that a countable product of compact spaces is a compact space: all we need to do is to verify the reverse induction step (\ding{70}) and this verification repeats the proof of compactness of the product of two compact spaces.

Note that since a countable product of Lindel\"{o}ff spaces is not always Lindel\"{o}ff, we cannot replace compactness by the Lindel\"{o}ff property in the above example of application. This is because the reverse induction step in this case is not valid.
But the principle of Reverse induction for Lindel\"{o}ff property \emph{is} valid (and it may be useful in the case of product of particular spaces for which the reverse induction step can be proved).

A similar situation occurs in the next section when we prove the D property of the countable power of the Sorgenfrey line. 
The principle of Reverse induction for D property (Theorem~\ref{teo.rev.ind}) is valid, the reverse induction step cannot be proved in the case of product of arbitrary D-spaces (because the D property does not preserved by products~\cite{D}), but it can be proved in the case of product of Sorgenfrey lines (see the proof of Theorem~\ref{teo.S.D}).

%
%
%
%
%
%
%
%
%

\section{D property of the countable power of the Sorgenfrey line}
\label{sect.S.is.D}

\begin{nota}\mbox{ }

  \begin{itemize}
  \item [\ding{46}\:]
    $\SS\:\coloneq\:\nos\,$the Sorgenfrey line.
  \item [\ding{46}\:]
    $[{a},{b})_\SS\:\coloneq\:\nos\,$ the half-interval $[{a},{b})$ as a subspace of $\SS$.
  \item [\ding{46}\:]
    $\SS^\omega\:\coloneq\:\prod_{{i}\in\omega}\SS_{i}\nos\,$,\quad
    where $\SS_{i}\:\coloneq\:[0,1)_\SS$ for all ${i}\in\omega$.
  \end{itemize}
\end{nota}

The spaces $[0,1)_\SS$ and $\SS$ can be written as a disjoint union of $\aleph_0$ many closed-and-open sets of the form $[{a},{b})_\SS$, so we have

\begin{rema}\label{rem.[0,1)}
$[0,1)_\SS$ is homeomorphic to $\SS$.\hfill$\qed$%
\end{rema}

\begin{nota}
Suppose that ${p}=\ll{p}_{i}\rr_{{i}\in\omega}\in\SS^\omega$, ${U}$ is a neighbourhood of ${p}$, ${m}\in\omega$, and $\varphi$ is an NA for $\SS^\omega$.
Then:

  \begin{itemize}
  \item [\ding{46}\:]

    ${U}\nos\,$ is a neighbourhood of ${p}$ of \emph{rank} ${m}\iff$

    ${U}\:\supseteq\:\big\{{q}=\ll{q}_{i}\rr_{{i}\in\omega}\in\SS^\omega: {q}_{i}\in[{p}_{i},{p}_{i}+2^{-m})_{\SS}
    \text{ for all }{i}\in{m}\big\}\nos\,$.
  \item [\ding{46}\:]
    $\mathsf{rank}({p},\varphi)\:\coloneq\:
    \mathsf{min}\{{m}\in\omega:\varphi({p})\text{ is a neighbourhood of }{p} \text{ of rank } {m}\}\nos\,$.
  \end{itemize}
\end{nota}

\noindent Note that each neighbourhood of a point in $\SS^\omega$ has some finite rank,
so $\mathsf{rank}({p},\varphi)$ is correctly defined.

The notion of rank resembles the notion of ${f}\nos$-index, which Peter de Caux used in~\cite{deCaux} to prove that the finite powers of the Sorgenfrey line are hereditarily D. We also use several other ideas from this paper.

\begin{teor}\label{teo.S.D}
  The countable power of the Sorgenfrey line is a D-space.
\end{teor}

The most difficult part of the proof of this theorem is Lemma~\ref{lem.F}, which is proved in Section~\ref{proof.of.lem.F}. The reverse induction hypothesis is used in the proof of Lemma~\ref{lem.G} in the same section.


\begin{proof}
It follows from Remark~\ref{rem.[0,1)} that the countable power of the Sorgenfrey line is homeomorphic to $\SS^\omega$.
Suppose that $\varphi$ is NA for $\SS^\omega$, ${n}\in\omega$, and ${P}$ is a closed $\nos{n}$-broom on $\SS^\omega$.
Using Theorem~\ref{teo.rev.ind} it is enough to show that $\varphi$ is correct on ${P}$
whenever $\mathsf{RIH}(\SS^\omega,\varphi,{n},{P})$ holds.

We shall construct a sequence $\ll{D}_{m}\rr_{{m}\in\omega}$ of subsets of $\SS^\omega$
such that, for all ${m}\in\omega$,

  \begin{itemize}
  \item [1.\,]
    ${D}_{m}\nos\,$ is a closed discrete set in $\SS^\omega$,
  \item [2.\,]
    ${D}_{m}\subseteq{P}\setminus\bigcup_{{i}<{m}}\varphi\ll{D}_{i}\rr\nos\,$,\quad and
  \item [3.\,]
    $\bigcup_{{i}\leqslant{m}}\varphi\ll{D}_{i}\rr
    \:\supseteq\:\{{p}\in{P}:\mathsf{rank}({p},\varphi)\leqslant{m}\}\nos\,$.
  \end{itemize}

First we show that (1)--(3) imply that $\varphi$ is  correct on ${P}$,
and after that build the a sequence $\ll{D}_{m}\rr_{{m}\in\omega}$
using $\mathsf{RIH}(\SS^\omega,\varphi,{n},{P})$.
Put ${D}\coloneq\bigcup_{{m}\in\omega}{D}_{m}$.
By (2), ${D}\subseteq{P}$,
and it follows from (3) that $\varphi\ll{D}\rr\supseteq{P}$.
We are to show that ${D}$ is a closed discrete set in the subspace ${P}$.

Suppose on the contrary that ${D}$ has an accumulation point ${r}$ in the subspace ${P}$;
then ${r}$ is an accumulation point of ${D}$ (in $\SS^\omega$).
Let ${m}'\coloneq\mathsf{rank}({r},\varphi)$
and ${U}\coloneq\bigcup_{{i}\leqslant{m}'}\varphi\ll{D}_{i}\rr$.
By (3), ${U}$ is an open neighbourhood of ${r}$ in $\SS^\omega$,
so ${r}$ is an accumulation point of the set ${U}\cap{D}$.
It follows from (2) that
${U}\cap\bigcup_{{i}>{m}'}{D}_{i}=\varnothing$,
so ${U}\cap{D}={U}\cap\bigcup_{{i}\leqslant{m}'}{D}_{i}$,
hence ${r}$ is an accumulation point of the set $\bigcup_{{i}\leqslant{m}'}{D}_{i}$.
But it follows from (1) that
$\bigcup_{{i}\leqslant{m}'}{D}_{i}$ is a closed discrete set in $\SS^\omega$, a contradiction.

To construct the sequence $\ll{D}_{m}\rr_{{m}\in\omega}$ we need the following proposition, which we prove in the next section.
\medskip

\begin{prop}\label{prop.main}
Suppose that $\varphi$ is an NA for $\SS^\omega$, ${n}\in\omega$, and ${P}$ is a closed $\nos{n}$-broom on $\SS^\omega$. Assume that $\mathsf{RIH}(\SS^\omega,\varphi,{n},{P})$ holds.
Then for every ${m}\in\omega$, every ${A}\subseteq{m}$, and every closed $({m}\!\setminus\!{A}\nos)$-broom ${R}\subseteq{P}\!$ on $\hspace{0.5pt}\SS^\omega$,
there exists a closed discrete set ${D}$ in $\SS^\omega$ such that

  \begin{itemize}
  \item [\ding{40}\,]
    ${D}\subseteq{R}\quad\nos\,$ and
  \item [\ding{40}\,]
    $\varphi\ll{D}\rr\supseteq
    \{{r}\in{R}:\mathsf{rank}({r},\varphi)={m}\}\nos\,$.
  \end{itemize}
\end{prop}
\medskip

We build $\ll{D}_{m}\rr_{{m}\in\omega}$ by recursion on ${m}$.
Assume that all ${D}_{i}$ with ${i}<{m}$ have been chosen.
Put ${R}\coloneq{P}\setminus\bigcup_{{i}<{m}}\varphi\ll{D}_{i}\rr$
and ${A}\coloneq{m}\setminus{n}$.
Since ${m}\setminus{A}={m}\setminus({m}\setminus{n})\subseteq{n}$
and ${R}\subseteq{P}$,
it follows from Remark~\ref{rem.brooms}(b,c)
that ${R}$ is a closed $\nos({m}\!\setminus\!{A})$-broom on $\SS^\omega$.
Then by Proposition~\ref{prop.main} there exists a closed discrete set ${D}_{m}$ in $\SS^\omega$ such that

  \begin{itemize}
  \item [\ding{40}\,]
    ${D}_{m}\subseteq{R}\quad\nos\,$ and
  \item [\ding{40}\,]
    $\varphi\ll{D}_{m}\rr\supseteq
    \{{r}\in{R}:\mathsf{rank}({r},\varphi)={m}\}\nos\,$,
  \end{itemize}
from which conditions (1)--(3) easily follow.
\end{proof}

\section{Proof of Proposition~\ref{prop.main}}


In the following notation $[a,b]$ denotes the closed interval in the real line with endpoints ${a}$ and ${b}$.

\begin{nota}\label{nota.cube}
Suppose that ${m}\in\omega$ and $\varnothing\neq{A}\subseteq{m}$.
Then

  \begin{itemize}
  \item [\ding{46}\:]
    ${\ZZ}\nos\,$ is an ${A},{m}\nos$-\emph{cube}$\iff$

    ${\ZZ}=\prod_{{i}\in{A}}[{f}({i}),{f}({i})+2^{-m})_\SS\nos\,$
    for some function ${f}\colon{A}\to [0,1-2^{-{m}}]$.
  \end{itemize}
\end{nota}

Note that many of the following notations depend on ${m}$, ${A}$, and $\ZZ$, but we write, for example, $\edge({T})$ instead of $\edge_{{m},{A},{\ZZ}}({T})$ for better readability.

\begin{nota}
Suppose that ${m}\in\omega$, $\varnothing\neq{A}\subseteq{m}$, and ${\ZZ}$ is an ${A},{m}\nos$-cube. Let ${x}=\ll{x}_{i}\rr_{{i}\in{A}}\in{\ZZ}$,  ${y}=\ll{y}_{i}\rr_{{i}\in{A}}\in{\ZZ}$, and ${B}\subseteq{A}$.
Then:

  \begin{itemize}
  \item [\ding{46}\:]
    ${x}\sqsubseteq{y}\iff\forall{i}\in{A}\:[{x}_{i}\leqslant{y}_{i}]\nos\,$.
  \item [\ding{46}\:]
    ${x}\sqsubseteq_{B}{y}\iff
    \forall{i}\in{B}\:[{x}_{i}={y}_{i}]\nos\ $ and
    $\ \forall{j}\in{A}\setminus{B}\:[{x}_{j}<{y}_{j}]$.
  \end{itemize}
\end{nota}

\begin{rema}\label{rem.sqsubset}
Suppose that ${m}\in\omega$, $\varnothing\neq{A}\subseteq{m}$, and ${\ZZ}$ is an ${A},{m}\nos$-cube. Let
${x},{y},{z}\in{\ZZ}$, ${Y}\subseteq{\ZZ}$, and ${B},{C}\subseteq{A}$.
Then:

  \begin{itemize}
  \item [\textup{a.}\,]
    ${x}\sqsubseteq{y}{\quad{\longleftrightarrow}\quad}
    \exists{D}\subseteq{A}\:[{x}\sqsubseteq_{D}{y}]$.
  \item [\textup{b.}\,]
    If ${x}\sqsubseteq_{B}{y}$ and ${y}\sqsubseteq_{C}{z}$, then
    ${x}\sqsubseteq_{{B}\cap{C}}{z}$.
  \item [\textup{c.}\,]
    If ${x}\sqsubseteq_{B}{y}$ and ${B}\neq{A}$, then
    ${x}\neq{y}$.
  \item [\textup{d.}\,]
    If ${y}\sqsubseteq_\varnothing{x}$ and ${y}\in\mathsf{closure}({Y},{\ZZ})$,

    then there is ${z}\in{Y}$ such that ${z}\sqsubseteq_\varnothing{x}$.\hfill$\qed$%
  \end{itemize}
\end{rema}

\begin{nota}
Suppose that ${m}\in\omega$, $\varnothing\neq{A}\subseteq{m}$, and ${\ZZ}$ is an ${A},{m}\nos$-cube. Let ${p}\in\SS^\omega$ and ${X}\subseteq{\ZZ}$.
Then:

  \begin{itemize}
  \item [\ding{46}\:]
    $\dot{p}\:\coloneq\:{p}\res{A}\nos\,$.
  \item [\ding{46}\:]
    $\breve{X}\:\coloneq\:
    \{{p}\in\SS^\omega:\dot{p}\in{X}\}\nos\,$.
  \end{itemize}
\end{nota}

\begin{rema}\label{rem.hat}
Suppose that ${m}\in\omega$, $\varnothing\neq{A}\subseteq{m}$, and ${\ZZ}$ is an ${A},{m}\nos$-cube. Let ${p}\in\SS^\omega$ and ${X},{Y}\subseteq{\ZZ}$.
Then:

  \begin{itemize}
  \item [\textup{a.}\,]
    $\dot{p}\in{X}\quad{\longleftrightarrow}\quad
    {p}\in\breve{X}\nos\,$.
  \item [\textup{b.}\,]
    ${X}\subseteq{Y}\quad{\longleftrightarrow}\quad
    \breve{X}\subseteq\breve{Y}\nos\,$.
  \item [\textup{c.}\,]
    $\breve\ZZ\nos\,$ is a closed set in $\SS^\omega$.\hfill$\qed$%
\end{itemize}
\end{rema}
\noindent
We shall use item (a) of the above remark very often, sometimes without mention.

\begin{nota}
Suppose that ${m}\in\omega$, $\varnothing\neq{A}\subseteq{m}$, and ${\ZZ}$ is an ${A},{m}\nos$-cube. Let ${x}\in{\ZZ}$ and ${T}\subseteq\breve\ZZ$.
Then:

  \begin{itemize}
  \item [\ding{46}\:]
    $\cone({x})\:\coloneq\:\{{y}\in{\ZZ}:{x}\sqsubseteq_\varnothing{y}\}\nos\,$.
  \item [\ding{46}\:]
    $\crown({T})\:\coloneq\:\bigcup\{\cone(\dot{t}):{t}\in{T}\}\nos\,$.
  \item [\ding{46}\:]
    $\edge({T})\:\coloneq\:\nos\,$ the boundary of $\crown({T})$ in ${\ZZ}$.
  \item [\ding{46}\:]
    $\bottom({T})\:\coloneq\:
    \mathsf{min}_{\sqsubseteq}\big(\edge({T})\big)\:=\:
    \big\{{x}\in\edge({T}):
    \neg\exists{y}\in\edge({T})\:[{x}\neq{y}\sqsubseteq{x}]\big\}\nos\,$.
  \end{itemize}
\end{nota}

\begin{rema}\label{rem.psi}
Suppose that ${m}\in\omega$, $\varnothing\neq{A}\subseteq{m}$, and ${\ZZ}$ is an ${A},{m}\nos$-cube. Let ${x}\in{\ZZ}$, ${p}\in\breve\ZZ$, and  ${T}\subseteq\breve\ZZ$.
Then:

  \begin{itemize}
  \item [\textup{a.}\,]
    If ${y}\in\crown({T})$,
    then $\exists{t}\in{T}\:[\dot{t}\sqsubseteq_\varnothing{y}]$.
  \item [\textup{b.}\,]
    If ${y}\in\crown({T})$ and ${y}\sqsubseteq{x}$,
    then ${x}\in\crown({T})$.
  \item [\textup{c.}\,]
    $\edge({T})\:=\:
    \mathsf{closure}\big(\crown({T}),{\ZZ}\big)\setminus\crown({T})\nos\,$.
  \item [\textup{d.}\,]
    $\hedge({T})\:\supseteq\:
    \mathsf{closure}\big(\hcrown({T}),\SS^\omega\big)\setminus\hcrown({T})\nos\,$.
  \item [\textup{e.}\,]
    If ${q}\in{T}$,
    then $\dot{q}\in\mathsf{closure}\big(\crown({T}),{\ZZ}\big)$.
  \item [\textup{f.}\,]
    If ${T}\cap\hedge({T})=\varnothing$,
    then ${T}\subseteq\hcrown({T})$.
  \item [\textup{g.}\,]
    If ${q}\in{T}$ and $\dot{q}\sqsubseteq{y}\in\edge({T})$,
    then $\dot{q}\in\edge({T})$.
  \item [\textup{h.}\,]
    If ${q}\in\mathsf{closure}({T},\SS^\omega)$,
    then $\dot{q}\in\mathsf{closure}\big(\crown({T}),{\ZZ}\big)$.
  \item [\textup{i.}\,]
    If ${y}\neq\dot{p}\sqsubseteq{y}\in\bottom({T})$,
    then ${p}\notin\mathsf{closure}({T},\SS^\omega)$.
  \item [\textup{j.}\,]
    If $\varphi$ is an NA for $\SS^\omega$, ${q}\in{Q}\subseteq\breve\ZZ$, ${Q}$ is an $({m}\!\setminus\!{A}\nos)$-broom on $\SS^\omega$, and $\mathsf{rank}({q},\varphi)={m}$,

    then $\varphi({q})\supseteq{Q}\cap\hcone(\dot{q})$.
\end{itemize}
\end{rema}

\begin{proof}
Part (c) follows from the fact that $\crown({T})$ is open in $\ZZ$;
(d) follows from (c) and Remark~\ref{rem.hat}(a);
(f) follows from (c), (e), and Remark~\ref{rem.hat}(a);
(g) follows from (e), (c), and (b);
(h) follows from (e);
(i) follows from (h), (c), (b), and the definition of $\bottom({T})$.
\end{proof}

We shall prove Proposition~\ref{prop.main} by induction on the cardinality of the set ${A}$. In the following notation the formula $\text{\ding{74}}(\varphi,{n},{P},{l})$ expresses the statement of Proposition~\ref{prop.main} for ${A}$ of cardinality ${l}$; the formula $\mathsf{IH}(\varphi,{n},{P},{l})$ expresses the Induction Hypotheses.
Also $[{m}]^{l}$ denotes the set of subsets of ${m}$ of cardinality ${l}$.


\begin{nota}\label{nota.IH}
Suppose that $\varphi$ is an NA for $\SS^\omega$, ${n}\in\omega$, ${P}$ is a closed $\nos{n}$-broom on $\SS^\omega$, and ${l}\in\omega$. Then:

  \begin{itemize}
  \item [\ding{46}\:]
    $\text{\ding{74}}(\varphi,{n},{P},{l})\iff\nos\,$

    For every ${m}\in\omega$, every ${A}\in[{m}]^{{l}}$, and every closed $({m}\!\setminus\!{A}\nos)$-broom ${R}\subseteq{P}$ on $\SS^\omega$,

    there exists a closed discrete set ${D}$ in $\SS^\omega$ such that

    \begin{itemize}
    \item [\ding{40}\,]
      ${D}\subseteq{R}\quad\nos\,$ and
    \item [\ding{40}\,]
      $\varphi\ll{D}\rr\supseteq
      \{{r}\in{R}:\mathsf{rank}({r},\varphi)={m}\}\nos\,$.
    \end{itemize}
  \item [\ding{46}\:]
    $\mathsf{IH}(\varphi,{n},{P},{l})\iff
    \forall{k}<{l}\ \text{\ding{74}}(\varphi,{n},{P},{k})\nos\,$.
  \end{itemize}
\end{nota}

\begin{proof}[\textup{\textbf{Proof of Proposition~\ref{prop.main}}}]
Suppose that $\varphi$ is an NA for $\SS^\omega$, ${n}\in\omega$, and ${P}$ is a closed $\nos{n}$-broom on $\SS^\omega$. Assume that $\mathsf{RIH}(\SS^\omega,\varphi,{n},{P})$ holds. We shall prove $\text{\ding{74}}(\varphi,{n},{P},{l})$ for all ${l}\in\omega$, from which the conclusion of Proposition~\ref{prop.main} follows.

We proceed by induction on ${l}$:
for each ${l}\in\omega$, we show that $\text{\ding{74}}(\varphi,{n},{P},{l})$
follows from $\mathsf{IH}(\varphi,{n},{P},{l})$.
Suppose that ${m}\in\omega$, ${A}\in[{m}]^{{l}}$, and ${R}\subseteq{P}$ is a closed $({m}\!\setminus\!{A}\nos)$-broom on $\SS^\omega$.
We must find a closed discrete set ${D}$ in $\SS^\omega$ such that

\begin{itemize}
\item [\ding{40}\,]
  ${D}\subseteq{R}\quad\nos\,$ and

\item [\ding{40}\,]
  $\varphi\ll{D}\rr\supseteq
  \{{r}\in{R}:\mathsf{rank}({r},\varphi)={m}\}\nos\,$.
\end{itemize}

\emph{Case 1.} ${l}=0$.
Then ${A}$ is empty, so ${R}$ is an ${m}\nos$-broom on $\SS^\omega$. If the set $\{{r}\in{R}:\mathsf{rank}({r},\varphi)={m}\}$ is empty, then ${D}\coloneq\varnothing$ satisfies the required conditions.
Otherwise fix ${r}\in{R}$ such that $\mathsf{rank}({r},\varphi)={m}$.
Then $\varphi({r})\supseteq{R}$, therefore ${D}\coloneq\{r\}\subseteq{R}$ works.

\emph{Case 2.} ${l}\geqslant1$.
\emph{Case 2.1.} ${A}\cap{n}\neq\varnothing$.
Since ${R}\subseteq{P}$ and ${P}$ is an ${n}\nos$-broom on $\SS^\omega$,
${R}$ is also an ${n}\nos$-broom by Remark~\ref{rem.brooms}(c),
so ${R}$ is an $\big(({m}\!\setminus\!{A})\cup{n}\big)\nos$-broom by Remark~\ref{rem.brooms}(d).
Put ${B}\coloneq{A}\setminus{n}$.
We have ${m}\setminus{B}={m}\setminus({A}\!\setminus\!{n})\subseteq
({m}\!\setminus\!{A})\cup{n}$,
hence ${R}$ is an $({m}\!\setminus\!{B})\nos$-broom by~Remark~\ref{rem.brooms}(b).
We have $|{B}|<|{A}|={l}$ because of (Case 2.1), so, by $\mathsf{IH}(\varphi,{n},{P},{l})$,
$\text{\ding{74}}(\varphi,{n},{P},|{B}|)$ holds.
Applying it to ${m}\in\omega$, ${B}\in[{m}]^{|{B}|}$, and the closed $({m}\!\setminus\!{B}\nos)$-broom ${R}\subseteq{P}$,
we receive the desired set ${D}$.

\emph{Case 2.2.} ${A}\cap{n}=\varnothing$.
The space $\SS^\omega$ can be covered by a family
$\{\breve\ZZ_{j}:{j}\in{J}\}$, where each ${\ZZ}_{j}$ is an ${A},{m}\nos$-cube
and ${J}$ is finite.
For each ${j}\in{J}$, put
$$
{R}_{j}\coloneq{R}\cap\breve\ZZ_{j}\qquad \text{and}\qquad {T}_{j}\coloneq\{{r}\in{R}_{j}:\mathsf{rank}({r},\varphi)={m}\}.
$$
We have

  \begin{itemize}
  \item [1.\,]
    ${R}=\bigcup_{{j}\in{J}}{R}_{j}\nos\,$\quad and
  \item [2.\,]
    ${R}_{j}\nos\,$ is closed in $\SS^\omega\quad$ for all ${j}\in{J}$
  \end{itemize}
by Remark~\ref{rem.hat}(c).

To proceed, we need a lemma, which we prove in the next section.
\medskip

\begin{lemm}\label{lemm.I}
Suppose that $\varphi$ is an NA for $\SS^\omega$, ${n}\in\omega$, ${P}$ is a closed $\nos{n}$-broom on $\SS^\omega$, and ${l}\in\omega\setminus\{0\}$.
Assume that $\mathsf{RIH}(\SS^\omega,\varphi,{n},{P})$ and $\mathsf{IH}(\varphi,{n},{P},{l})$ hold.
Let ${m}\in\omega$, ${m}>{n}$, ${A}\in[{m}]^l$, ${\ZZ}$ be an ${A},{m}\nos$-cube, and ${T}\subseteq{P}\cap\breve\ZZ$ be an $({m}\!\setminus\!{A}\nos)$-broom on $\SS^\omega$ such that $\mathsf{rank}({t},\varphi)={m}$ for all ${t}\in{T}$.

Then there are ${Q}\subseteq{T}$ and a closed discrete set ${D}$ in $\SS^\omega$
such that

  \begin{itemize}
  \item [\textup{a.}\,]
    ${D}\subseteq\mathsf{closure}({T},\SS^\omega)\nos\,$,
  \item [\textup{b.}\,]
    $\varphi\ll{D}\rr\supseteq{T}\setminus{Q}\nos\,$,
  \item [\textup{c.}\,]
    ${Q}\subseteq\hcrown({Q})\nos\,$,\quad and
  \item [\textup{d.}\,]
    $\mathsf{closure}({Q},\SS^\omega)\cap
    \hbottom({Q})=\varnothing\nos\,$.
  \end{itemize}
\end{lemm}
\medskip

We have ${A}\subseteq{m}$, ${A}\neq\varnothing$ by (Case 2),
and ${A}\cap{n}=\varnothing$ by (Case 2.2);
therefore ${m}>{n}$.
For each ${j}\in{J}$, ${T}_{j}\subseteq{R}_{j}={R}\cap\breve\ZZ_{j}$, ${R}\subseteq{P}$,
and ${R}$ is an $({m}\!\setminus\!{A})\nos$-broom on $\SS^\omega$,
so ${T}_{j}\subseteq{P}\cap\breve\ZZ_{j}$ and
${T}_{j}$ is an $({m}\!\setminus\!{A})\nos$-broom by Remark~\ref{rem.brooms}(c).
Also, for each ${j}\in{J}$, we have $\mathsf{rank}({t},\varphi)={m}$ for all ${t}\in{T}_{j}$ by the choice of ${T}_{j}.$
So we may apply Lemma~\ref{lemm.I} to $\ZZ=\ZZ_{j}$ and ${T}={T}_{j}$.
By this lemma, for each ${j}\in{J}$, there are
${Q}_{j}\subseteq{T}_{j}$ and a closed discrete set ${D}_{j}$ in $\SS^\omega$
such that

  \begin{itemize}
  \item [A.\,]
    ${D}_{j}\subseteq\mathsf{closure}({T}_{j},\SS^\omega)\nos\,$,
  \item [B.\,]
    $\varphi\ll{D}_{j}\rr\supseteq{T}_{j}\setminus{Q}_{j}\nos\,$,
  \item [C.\,]
    ${Q}_{j}\subseteq\hcrown({Q}_{j})\nos\,$,\quad and
  \item [D.\,]
    $\mathsf{closure}({Q}_{j},\SS^\omega)\cap
    \hbottom({Q}_{j})=\varnothing\nos\,$.
  \end{itemize}

Now we need a second lemma, which we also prove in the next section.
Note that we do not use neither $\mathsf{RIH}(\SS^\omega,\varphi,{n},{P})$ nor $\mathsf{IH}(\varphi,{n},{P},{l})$ in the assumption of this lemma.

\medskip

\begin{lemm}\label{lemm.II}
Let $\varphi$ be an NA for $\SS^\omega$,
${m}\in\omega$, $\varnothing\neq{A}\subseteq{m}$, ${\ZZ}$ be an ${A},{m}\nos$-cube,
and ${Q}\subseteq\breve\ZZ$ be an $({m}\!\setminus\!{A}\nos)$-broom on $\SS^\omega$ such that $\mathsf{rank}({q},\varphi)={m}$ for all ${q}\in{Q}$.
Suppose that

  \begin{itemize}
  \item [\textup{c.}\,]
    ${Q}\subseteq\hcrown({Q})\quad\nos\,$ and
  \item [\textup{d.}\,]
    $\mathsf{closure}({Q},\SS^\omega)\cap
    \hbottom({Q})=\varnothing\nos\,$.
  \end{itemize}

Then there is a closed discrete set ${E}$ in $\SS^\omega$ such that

  \begin{itemize}
  \item [\textup{e.}\,]
    ${E}\subseteq{Q}\quad\nos\,$ and
  \item [\textup{f.}\,]
    $\varphi\ll{E}\rr\supseteq{Q}\nos\,$.
  \end{itemize}
\end{lemm}
\medskip

Using (C) and (D), we may apply Lemma~\ref{lemm.II} to $\ZZ=\ZZ_{j}$ and ${Q}={Q}_{j}$
because ${A}\neq\varnothing$ by (Case 2) and, for each ${j}\in{J}$,
${Q}_{j}\subseteq{T}_{j}\subseteq\breve\ZZ_{j}$
and ${T}_{j}$ is an $({m}\!\setminus\!{A}\nos)$-broom on $\SS^\omega$ such that $\mathsf{rank}({t},\varphi)={m}$ for all ${t}\in{T}_{j}$.
By this lemma, for each ${j}\in{J}$,
there is a closed discrete set ${E}_{j}$ in $\SS^\omega$ such that

  \begin{itemize}
  \item [E.\,]
    ${E}_{j}\subseteq{Q}_{j}\quad\nos\,$ and
  \item [F.\,]
    $\varphi\ll{E}_{j}\rr\supseteq{Q}_{j}\nos\,$.
  \end{itemize}

Put ${D}\coloneq\bigcup_{{j}\in{J}}({D}_{j}\cup{E}_{j})$.
Then ${D}$ is a closed discrete set in $\SS^\omega$ because all ${D}_{j}$ and ${E}_{j}$ are such and ${J}$ is finite.
Using (A), (E), ${Q}_{j}\subseteq{T}_{j}$, (2), and (1), we have
$$
{D}\:\subseteq\:
\bigcup_{{j}\in{J}}\big(\mathsf{closure}({T}_{j},\SS^\omega)\cup{Q}_{j}\big)
\:\subseteq\:
\bigcup_{{j}\in{J}}\mathsf{closure}({T}_{j},\SS^\omega)
\:\subseteq\:
\bigcup_{{j}\in{J}}\mathsf{closure}({R}_{j},\SS^\omega)
\:=\:
\bigcup_{{j}\in{J}}{R}_{j}\:=\:{R},
$$
so ${D}\subseteq{R}$.
At last, using (B), (F), and (1), we get
$$
\varphi\ll{D}\rr\:=\:
\bigcup_{{j}\in{J}}\varphi\ll{D}_{j}\cup{E}_{j}\rr
\:\supseteq\:
\bigcup_{{j}\in{J}}\big(({T}_{j}\setminus{Q}_{j})\cup{Q}_{j}\big)
\:\supseteq\:
\bigcup_{{j}\in{J}}{T}_{j}
\:\supseteq\:
\big\{{r}\in{\textstyle\bigcup_{{j}\in{J}}}{R}_{j}={R}:\mathsf{rank}({r},\varphi)={m}\big\}.
$$%
\end{proof}%

\section{Proofs of Lemmas~\ref{lemm.I} and~\ref{lemm.II}}

\begin{nota}
Suppose that ${m}\in\omega$, $\varnothing\neq{A}\subseteq{m}$, and ${\ZZ}$ is an ${A},{m}\nos$-cube. Let ${B}\subseteq{A}$, ${y}\in{\ZZ}$, and
${z}_{i}\in{\ZZ}$ for all ${i}\in\omega$.
Then

  \begin{itemize}
  \item [\ding{46}\:]
    ${z}_{i}\xrightarrow[{i}\to\infty]{{B}\downarrow}{y}\,
    \iff\nos\,$
    the sequence $\ll{z}_{i}\rr_{{i}\in\omega}\nos\,$ converges to ${y}$ in ${\ZZ}$ and
    ${y}\sqsubseteq_{B}{z}_{{i}+1}\sqsubseteq_{B}{z}_{i}$ for all ${i}\in\omega$.
  \end{itemize}
\end{nota}

\begin{rema}\label{rem.1.25}
Suppose that ${m}\in\omega$, $\varnothing\neq{A}\subseteq{m}$, and ${\ZZ}$ is an ${A},{m}\nos$-cube. Let ${x},{y}\in{\ZZ}$, ${C},{D}\subseteq{A}$,
${x}\sqsubseteq_{{C}\cap{D}}{y}$, and $\ll{z}_{j}\rr_{{j}\in\omega}$ be a sequence in ${\ZZ}$ such that ${z}_{j}\xrightarrow[{j}\to\infty]{{D}\downarrow}{x}$.

Then there is ${j}'\in\omega$ such that
${z}_{{j}'}\sqsubseteq_{{C}\cap{D}}{y}$.\hfill$\qed$%
\end{rema}

If a sequence $\ll{z}_{i}\rr_{{i}\in\omega}\nos\,$ in $\SS$ converges to a point ${y}$,
then it has a subsequence $\ll{z}_{{i}_{k}}\rr_{{k}\in\omega}$ such that
either each ${z}_{{i}_{k}}$ equals ${y}$
or ${y}<{z}_{{i}_{{k}+1}}<{z}_{{i}_{k}}$ for all ${k}\in\omega$.
Repeating this argument $|{A}|$ times, we receive the following generalization:

\begin{rema}\label{lemm.C-lim}
Suppose that ${m}\in\omega$, $\varnothing\neq{A}\subseteq{m}$, and ${\ZZ}$ is an ${A},{m}\nos$-cube. Let ${r}\in\SS^\omega$ and $\ll{p}_{i}\rr_{{i}\in\omega}$ be a sequence in $\breve\ZZ$ that converges to ${r}$ in $\SS^\omega$.

Then $\dot{r}\in{\ZZ}$ and there are ${B}\subseteq{A}$ and a subsequence $\ll{p}_{{i}_{k}}\rr_{{k}\in\omega}$ such that $\dot{p}_{{i}_{k}}\xrightarrow[{k}\to\infty]{{B}\downarrow}\dot{r}$.\hfill$\qed$%
\end{rema}

\begin{proof}[\textbf{Proof of Lemma~\ref{lemm.II}}]
Let $\varphi$ be an NA for $\SS^\omega$, ${m}\in\omega$, $\varnothing\neq{A}\subseteq{m}$, ${\ZZ}$ be an ${A},{m}\nos$-cube, and ${Q}\subseteq\breve\ZZ$ be an $({m}\!\setminus\!{A}\nos)$-broom on $\SS^\omega$ such that $\mathsf{rank}({q},\varphi)={m}$ for all ${q}\in{Q}$.
Suppose that

  \begin{itemize}
  \item [\textup{c.}\,]
    ${Q}\subseteq\hcrown({Q})\quad\nos\,$ and
  \item [\textup{d.}\,]
    $\mathsf{closure}({Q},\SS^\omega)\cap
    \hbottom({Q})=\varnothing\nos\,$.
  \end{itemize}
We must find a closed discrete set ${E}$ in $\SS^\omega$ such that

  \begin{itemize}
  \item [\textup{e.}\,]
    ${E}\subseteq{Q}\quad\nos\,$ and
  \item [\textup{f.}\,]
    $\varphi\ll{E}\rr\supseteq{Q}\nos\,$.
  \end{itemize}

If ${Q}$ is empty, then we can take ${E}\coloneq\varnothing$; so we consider the case when ${Q}\neq\varnothing$.
Since $\crown({Q})=\bigcup\{\cone(\dot{p}):{p}\in{Q}\}$ and each set $\cone(\dot{p})$ is open in the Euclidean topology on ${\ZZ}$, which is hereditarily Lindel\"{o}ff, there is a sequence $\ll{p}_{i}\rr_{{i}\in\omega}$ in ${Q}$ such that $\crown({Q})\subseteq\bigcup\{\cone(\dot{p}_{i}):{i}\in\omega\}$.
For ${x}=\ll{x}_{j}\rr_{{j}\in{A}}\in{\ZZ}$ and ${y}=\ll{y}_{j}\rr_{{j}\in{A}}\in{\ZZ}$,
let $\rho_\Sigma({x},{y})\:\coloneq\:\sum_{{j}\in{A}}\|{x}_{j}-{y}_{j}\|$,
where $\|a\|$ is the absolute value of a real number ${a}$.

We shall build a sequence $\ll{q}_{i}\rr_{{i}\in\omega}$ in ${Q}$ such that,
for all ${i},{j}\in\omega$,

  \begin{itemize}
  \item [1.\,]
    $\dot{q}_{i}\sqsubseteq_\varnothing\dot{p}_{i}\nos\,$,
  \item [2.\,]
    if ${j}<{i}$ and $\dot{q}_{i}\sqsubseteq\dot{q}_{j}$,
    then $\dot{q}_{i}\sqsubseteq_\varnothing\dot{q}_{j}$,\quad and
  \item [3.\,]
    if ${s}\in{Q}$ and $\dot{s}\sqsubseteq_\varnothing\dot{q}_{i}$, then $\rho_{\Sigma}(\dot{s},\dot{q}_{i})<{i}^{-1}$.
  \end{itemize}

Put ${E}\coloneq\{{q}_{i}:{i}\in\omega\}$. First we show that (1)--(3) imply that ${E}$
has the required properties, and after that build the sequence $\ll{q}_{i}\rr_{{i}\in\omega}$.
Condition (e) holds because all ${q}_{i}$ are in ${Q}$.
Since $\mathsf{rank}({q},\varphi)={m}$ for all ${q}\in{Q}$,
it follows from Remark~\ref{rem.psi}(j)
that $\varphi({q}_{i})\supseteq{Q}\cap\hcone(\dot{q}_{i})$
for all ${i}\in\omega$;
also (1) with Remark~\ref{rem.sqsubset}(b) imply $\cone(\dot{q}_{i})\supseteq\cone(\dot{p}_{i})$ for all ${i}\in\omega$.
Therefore, using Remark~\ref{rem.hat}(b), the choice of the sequence $\ll{p}_{i}\rr_{{i}\in\omega}$, Remark~\ref{rem.hat}(a), and (c), we can write
$$
\varphi\ll{E}\rr\:=\:
\bigcup_{{i}\in\omega}\varphi({q}_{i})
\:\supseteq\:
\bigcup_{{i}\in\omega}\big({Q}\cap\hcone(\dot{q}_{i})\big)
\:\supseteq\:
$$
$$
{Q}\cap\bigcup_{{i}\in\omega}\hcone(\dot{q}_{i})
\:\supseteq\:
{Q}\cap\bigcup_{{i}\in\omega}\hcone(\dot{p}_{i})
\:\supseteq\:
{Q}\cap\hcrown({Q})
\:\supseteq\:{Q}\cap{Q}\:=\:{Q},
$$
so (f) also holds.
We are to show that ${E}$ is a closed discrete set in $\SS^\omega$.

Suppose on the contrary that ${E}$  has an accumulation point ${r}$ in $\SS^\omega$.
Using a countable local base at the point ${r}$ in $\SS^\omega$, we can build a strictly increasing sequence $\ll{i}_{k}\rr_{{k}\in\omega}$ in $\omega$ such that $\{{q}_{{i}_{k}}:{k}\in\omega\}\subseteq{E}\setminus\{{r}\}$ and
the subsequence $\ll{q}_{{i}_{k}}\rr_{{k}\in\omega}$ converges to~${r}$.
Using Remark~\ref{lemm.C-lim} and passing to a subsequence,
we may assume that $\dot{q}_{{i}_{k}}\xrightarrow[{k}\to\omega]{{B}\downarrow}\dot{r}$
for some ${B}\subseteq{A}$.
Then it follows from (2) and Remark~\ref{rem.sqsubset}(a,b) that $\dot{q}_{{i}_{k}}\xrightarrow[{k}\to\omega]{\varnothing\downarrow}\dot{r}$.

Since ${E}\subseteq{Q}$, we have ${r}\in\mathsf{closure}({Q},\SS^\omega)$, so
(d) implies ${r}\notin\hbottom({Q})$.
We also have ${r}\notin\hcrown({Q})$:
If not, then
$\dot{r}\in\crown({Q})$ by Remark~\ref{rem.hat}(a),
so, by Remark~\ref{rem.psi}(a),
$\dot{s}\sqsubseteq_\varnothing\dot{r}$ for some ${s}\in{Q}$.
Then $\dot{s}\sqsubseteq_\varnothing\dot{r}\sqsubseteq_\varnothing\dot{q}_{{i}_{k}}$,
hence $\dot{s}\sqsubseteq_\varnothing\dot{q}_{{i}_{k}}$, for all ${k}\in\omega$.
From this (3) implies $\rho_{\Sigma}(\dot{s},\dot{q}_{{i}_{k}})<{{i}_{k}}^{-1}$
for all ${k}\in\omega$.
Since $\ll{i}_{k}\rr_{{k}\in\omega}$ strictly increases, this contradicts $\rho_{\Sigma}(\dot{s},\dot{q}_{{i}_{k}})
>\rho_{\Sigma}(\dot{s},\dot{r})>0$,
which also follows from $\dot{s}\sqsubseteq_\varnothing\dot{r}
\sqsubseteq_\varnothing\dot{q}_{{i}_{k}}$.
Using all above with (c) and Remark~\ref{rem.psi}(d) we get
$$
{r}\:\in\:\big(\mathsf{closure}({Q},\SS^\omega)\setminus\hcrown({Q})\big)
\setminus\hbottom({Q})
\:\subseteq\:
$$
$$
\Big(\mathsf{closure}\big(\hcrown({Q}),\SS^\omega\big)\setminus\hcrown({Q})\Big)
\setminus\hbottom({Q})
\:\subseteq\:
\hedge({Q})\setminus\hbottom({Q}).
$$
Then $\dot{r}\in\edge({Q})\setminus\bottom({Q})$ by Remakr~\ref{rem.hat}(a),
so, by definition of $\bottom({Q})$, there is ${y}\in\edge({Q})$ such that $\dot{r}\neq{y}\sqsubseteq\dot{r}$.
Then $\rho_{\Sigma}({y},\dot{r})>0$,
so there is ${{k}'}\in\omega$ such that
$$
{{i}_{{k}'}}^{-1}<\rho_\Sigma({y},\dot{r})/3
\quad\text{and}\quad
\rho_\Sigma(\dot{q}_{{i}_{{k}'}},\dot{r})<\rho_\Sigma({y},\dot{r})/3
$$
because $\ll{i}_{k}\rr_{{k}\in\omega}$ strictly increases and
$\dot{q}_{{i}_{k}}\xrightarrow[{k}\to\omega]{\varnothing\downarrow}\dot{r}$.
We have ${y}\sqsubseteq\dot{r}\sqsubseteq_\varnothing\dot{q}_{{i}_{{k}'}}$, so
${y}\sqsubseteq_\varnothing\dot{q}_{{i}_{{k}'}}$.
Let
$$
{Y}\coloneq\{{x}\in\crown({Q}):\rho_\Sigma({y},{x})<\rho_\Sigma({y},\dot{r})/3\}.
$$
Since ${y}\in\edge({Q})\subseteq\mathsf{closure}\big(\crown({Q}),{\ZZ}\big)$,
we have ${y}\in\mathsf{closure}({Y},{\ZZ})$,
so by Remark~\ref{rem.sqsubset}(d)
there is ${z}\in{Y}$ such that ${z}\sqsubseteq_\varnothing\dot{q}_{{i}_{{k}'}}$.
By Remark~\ref{rem.psi}(a), since ${z}\in{Y}\subseteq\crown({Q})$,
there is ${s}\in{Q}$ such that $\dot{s}\sqsubseteq_\varnothing{z}$.
Then $\dot{s}\sqsubseteq_\varnothing\dot{q}_{{i}_{{k}'}}$,
so $\rho_\Sigma(\dot{s},\dot{q}_{{i}_{{k}'}})<
{{i}_{{k}'}}^{-1}<\rho_\Sigma({y},\dot{r})/3$ by (3) and the choice of ${{k}'}$.
Also $\dot{s}\sqsubseteq_\varnothing{z}\sqsubseteq_\varnothing\dot{q}_{{i}_{{k}'}}$
imply $\rho_\Sigma({z},\dot{q}_{{i}_{{k}'}})<
\rho_\Sigma(\dot{s},\dot{q}_{{i}_{{k}'}})<\rho_\Sigma({y},\dot{r})/3$;
and $\rho_\Sigma({y},{z})<\rho_\Sigma({y},\dot{r})/3$ because ${z}\in{Y}$.
Now, using the triangle inequality for the metric $\rho_\Sigma$, the above inequalities,
and the choice of ${k}'$, we have
$$
\rho_\Sigma({y},\dot{r})\:\leqslant\:
\rho_\Sigma({y},{z})+\rho_\Sigma({z},\dot{q}_{{i}_{{k}'}})
+\rho_\Sigma(\dot{q}_{{i}_{{k}'}},\dot{r})\:<\:
\rho_\Sigma({y},\dot{r})/3+\rho_\Sigma({y},\dot{r})/3+\rho_\Sigma({y},\dot{r})/3
\:=\:\rho_\Sigma({y},\dot{r}),
$$
a contradiction.

It remains to construct a sequence $\ll{q}_{i}\rr_{{i}\in\omega}$ in ${Q}$ that satisfies (1)--(3); we build it by recursion on~${i}$. Suppose that all ${q}_{j}$ with ${j}<{i}$ have been chosen.
It follows from (c) and Remark~\ref{rem.hat}(a) that $\dot{p}_{i}\in\crown({Q})$,
so by Remark~\ref{rem.psi}(a), there is ${t}\in{Q}$ such that $\dot{t}\sqsubseteq_\varnothing\dot{p}_{i}$.
Since ${Q}\subseteq\breve\ZZ$, $\dot{p}_{i}\in{\ZZ}$, and ${\ZZ}$ has a finite diameter in the metric $\rho_\Sigma$, among such ${t}$ there is one --- denote it by ${t}'$ --- with the additional property
$$
\mathsf{sup}\{\rho_{\Sigma}(\dot{s},\dot{p}_{i}):
{s}\in{Q}\text{ and }\dot{s}\sqsubseteq_\varnothing\dot{p}_{i}\}
-\rho_{\Sigma}(\dot{t}',\dot{p}_{i})
\:<\:{i}^{-1}.
$$
Repeating the argument by which we found ${t}\in{Q}$ at most ${i}+1$ times (when choosing a point ${u}_{{l}+1}\in{Q}$ in condition (4) below), it is straightforward to build, by recursion on ${l}$, a sequence $\ll{u}_{l}\rr_{{l}\leqslant{i}+1}$ in ${Q}$ such that
${u}_0={t}'\:(\in{Q})$ and, for all ${l}\leqslant{i}$,

  \begin{itemize}
  \item [4.\,]
    if there is ${j}<{i}$ such that ($\dot{u}_{l}\sqsubseteq\dot{q}_{j}$
    and $\dot{u}_{l}\not\sqsubseteq_\varnothing\dot{q}_{j}$),

    then $\dot{u}_{{l}+1}\sqsubseteq_\varnothing\dot{u}_{l}$;
  \item [5.\,]
    if there is no ${j}<{i}$ such that ($\dot{u}_{l}\sqsubseteq\dot{q}_{j}$
    and $\dot{u}_{l}\not\sqsubseteq_\varnothing\dot{q}_{j}$),

    then ${u}_{{l}+1}={u}_{l}$.
  \end{itemize}
Put ${q}_{i}\coloneq{u}_{{i}+1}$. Note that
$$
\dot{p}_{i}\sqsupseteq_\varnothing\dot{t}'=\dot{u}_0
\sqsupseteq_\varnothing\dots\sqsupseteq_\varnothing
\dot{u}_{k}=\dots=\dot{u}_{{i}+1}=\dot{q}_{i}
$$
for some ${k}\leqslant{i}+1$.
Now, condition (1) follows from Remark~\ref{rem.sqsubset}(b).
Also Remark~\ref{rem.sqsubset}(a,b) implies that
each ${j}<{i}$ meets the premise of condition (4) at most once
(i.e., for at most one ${l}$).
Since there are only ${i}$ different ${j}$ that are less than ${i}$,
at some step ${l}'\leqslant{i}$ of the recursion
no ${j}$ will meet the premise of (4);
that is,

  \begin{itemize}
  \item [6.\,]
    there is no ${j}<{i}$ such that ($\dot{u}_{{l}'}\sqsubseteq\dot{q}_{j}$
    and $\dot{u}_{{l}'}\not\sqsubseteq_\varnothing\dot{q}_{j}$).
  \end{itemize}
Hence, at step ${l}'$, the premise of (5) is satisfied,
so ${u}_{{l}'+1}={u}_{{l}'}$. Then, for the same reason,
$$
{u}_{{l}'+2}={u}_{{l}'+1},\qquad{u}_{{l}'+3}={u}_{{l}'+2},\qquad\dots\qquad,{u}_{{i}+1}={u}_{i}.
$$
It follows that $\dot{q}_{i}=\dot{u}_{{l}'}$, so, by (6),
there is no ${j}<{i}$ such that ($\dot{q}_{i}\sqsubseteq\dot{q}_{j}$ and $\dot{q}_{i}\not\sqsubseteq_\varnothing\dot{q}_{j}$).
This means that condition (2) is satisfied.

It remains to establish (3). Suppose that
${s}\in{Q}$ and $\dot{s}\sqsubseteq_\varnothing\dot{q}_{i}$.
Then $\dot{s}\sqsubseteq_\varnothing\dot{q}_{i}
\sqsubseteq\dot{t}'\sqsubseteq_\varnothing\dot{p}_{i}$,
so $\dot{s}\sqsubseteq_\varnothing\dot{p}_{i}$,
hence $\rho_{\Sigma}(\dot{s},\dot{p}_{i})
-\rho_{\Sigma}(\dot{t}',\dot{p}_{i})<{i}^{-1}$
by the choice of ${t}'$.
Also it follows from $\dot{s}\sqsubseteq_\varnothing\dot{q}_{i}
\sqsubseteq\dot{t}'\sqsubseteq_\varnothing\dot{p}_{i}$
that $\rho_{\Sigma}(\dot{s},\dot{p}_{i})=
\rho_{\Sigma}(\dot{s},\dot{t}')+\rho_{\Sigma}(\dot{t}',\dot{p}_{i})$
and $\rho_{\Sigma}(\dot{s},\dot{q}_{i})\leqslant\rho_{\Sigma}(\dot{s},\dot{t}')$.
Then
$$
\rho_{\Sigma}(\dot{s},\dot{q}_{i})\leqslant
\rho_{\Sigma}(\dot{s},\dot{t}')=
\rho_{\Sigma}(\dot{s},\dot{p}_{i})-\rho_{\Sigma}(\dot{t}',\dot{p}_{i})
<{i}^{-1}\,\text{,}
$$
so condition (3) is also satisfied.
\end{proof}

To prove Lemma~\ref{lemm.I} we need two lemmas, which we prove in the next section.
Note that the first of these lemmas uses only $\mathsf{IH}(\varphi,{n},{P},{l})$ as an assumption, and the second uses only $\mathsf{RIH}(\SS^\omega,\varphi,{n},{P})$.%
\smallskip

\begin{lemm}\label{lem.F}
Suppose that $\varphi$ is an NA for $\SS^\omega$, ${n}\in\omega$, ${P}$ is a closed $\nos{n}$-broom on $\SS^\omega$, and ${l}\in\omega\setminus\{0\}$.
Assume that $\mathsf{IH}(\varphi,{n},{P},{l})$ holds.
Let ${m}\in\omega$, ${A}\in[{m}]^l$, ${\ZZ}$ be an ${A},{m}\nos$-cube, and ${T}\subseteq{P}\cap\breve\ZZ$ be an $({m}\!\setminus\!{A}\nos)$-broom on $\SS^\omega$
such that $\mathsf{rank}({t},\varphi)={m}$ for all ${t}\in{T}$.

Then there is a closed discrete set ${F}$ in $\SS^\omega$ such that

  \begin{itemize}
  \item [\textup{g.}\,]
    ${F}\subseteq\mathsf{closure}
    \big({T}\cap\hedge({T}),\SS^\omega\big)\quad\nos\,$ and
  \item [\textup{h.}\,]
    $\varphi\ll{F}\rr\supseteq{T}\cap\hedge({T})\nos\,$.
  \end{itemize}
\end{lemm}
\smallskip

\begin{lemm}\label{lem.G}
Suppose that $\varphi$ is an NA for $\SS^\omega$, ${n}\in\omega$, and ${P}$ is a closed $\nos{n}$-broom on $\SS^\omega$.
Assume that $\mathsf{RIH}(\SS^\omega,\varphi,{n},{P})$ holds.
Let ${m}\in\omega$, ${m}>{n}$, $\varnothing\neq{A}\subseteq{m}$, ${\ZZ}$ be an ${A},{m}\nos$-cube, and ${T}\subseteq{P}\cap\breve\ZZ$ be an $({m}\!\setminus\!{A}\nos)$-broom on $\SS^\omega$.

Then there is a closed discrete set ${G}$ in $\SS^\omega$ such that

  \begin{itemize}
  \item [\textup{i.}\,]
    ${G}\subseteq
    \mathsf{closure}({T},\SS^\omega)\cap\hbottom({T})\quad\nos\,$ and
  \item [\textup{j.}\,]
    $\varphi\ll{G}\rr\supseteq
    \mathsf{closure}({T},\SS^\omega)\cap\hbottom({T})\nos\,$.
  \end{itemize}
\end{lemm}
\smallskip

\begin{proof}[\textbf{Proof of Lemma~\ref{lemm.I}}]
Suppose that $\varphi$ is an NA for $\SS^\omega$, ${n}\in\omega$, ${P}$ is a closed $\nos{n}$-broom on $\SS^\omega$, and ${l}\in\omega\setminus\{0\}$.
Assume that $\mathsf{RIH}(\SS^\omega,\varphi,{n},{P})$ and $\mathsf{IH}(\varphi,{n},{P},{l})$ hold.
Let ${m}\in\omega$, ${m}>{n}$, ${A}\in[{m}]^l$, ${\ZZ}$ be an ${A},{m}\nos$-cube, and ${T}\subseteq{P}\cap\breve\ZZ$ be an $({m}\!\setminus\!{A}\nos)$-broom on $\SS^\omega$ such that $\mathsf{rank}({t},\varphi)={m}$ for all ${t}\in{T}$.

\noindent
We must find ${Q}\subseteq{T}$ and a closed discrete set ${D}$ in $\SS^\omega$
such that

  \begin{itemize}
  \item [\textup{a.}\,]
    ${D}\subseteq\mathsf{closure}({T},\SS^\omega)\nos\,$,
  \item [\textup{b.}\,]
    $\varphi\ll{D}\rr\supseteq{T}\setminus{Q}\nos\,$,
  \item [\textup{c.}\,]
    ${Q}\subseteq\hcrown({Q})\nos\,$,\quad and
  \item [\textup{d.}\,]
    $\mathsf{closure}({Q},\SS^\omega)\cap
    \hbottom({Q})=\varnothing\nos\,$.
  \end{itemize}

Let $\delta$ be an ordinal of cardinality grater than the cardinality of ${T}$. Using Lemmas~\ref{lem.F} and~\ref{lem.G} with Remark~\ref{rem.brooms}(c), it is straightforward to build, by recursion on $\alpha$, three transfinite sequences $\ll{T}_\alpha\rr_{\alpha<\delta}$,
$\ll{F}_\alpha\rr_{\alpha<\delta}$, and $\ll{G}_\alpha\rr_{\alpha<\delta}$ such that, for all $\alpha<\delta$,

  \begin{itemize}
  \item [1.\,]
    ${T}_0={T}\nos\,$,
  \item [2.\,]
    ${T}_{\alpha+1}={T}_\alpha\setminus\varphi
    \ll{F}_\alpha\cup{G}_\alpha\rr\nos\,$,
  \item [3.\,]
    ${T}_\alpha=\bigcap_{\beta<\alpha}{T}_\beta\quad\nos\,$
    if $\alpha$ is a limit ordinal,
  \item [4.\,]
    ${F}_\alpha\nos\,$ and ${G}_\alpha$ are closed discrete sets in $\SS^\omega$,
  \item [G.\,]
    ${F}_\alpha\subseteq\mathsf{closure}
    \big({T}_\alpha\cap\hedge({T}_\alpha),\SS^\omega\big)\nos\,$,
  \item [H.\,]
    $\varphi\ll{F}_\alpha\rr\supseteq{T}_\alpha\cap
    \hedge({T}_\alpha)\nos\,$,
  \item [I.\,]
    ${G}_\alpha\subseteq
    \mathsf{closure}({T}_\alpha,\SS^\omega)\cap
    \hbottom({T}_\alpha)\nos\,$,\quad and
  \item [J.\,]
    $\varphi\ll{G}_\alpha\rr\supseteq
    \mathsf{closure}({T}_\alpha,\SS^\omega)\cap
    \hbottom({T}_\alpha)\nos\,$.
  \end{itemize}

Note that conditions (2) and (3) imply

  \begin{itemize}
  \item [5.\,]
    ${T}_\alpha\supseteq{T}_\beta\quad\nos\,$ for all $\alpha\leqslant\beta<\delta$.
  \end{itemize}
Also, (5), (1), (G), (I), and Remark~\ref{rem.hat}(c) imply

  \begin{itemize}
  \item [6.\,]
    ${T}_\alpha\subseteq\breve\ZZ\quad\nos\,$ for all $\alpha<\delta$;
  \item [7.\,]
    ${F}_\alpha\cup{G}_\alpha\subseteq\breve\ZZ\quad\nos\,$ for all $\alpha<\delta$.
  \end{itemize}

It follows from the choice of $\delta$, (1), and (5) that there is $\alpha<\delta$ such that ${T}_{\alpha+1}={T}_\alpha$.
Let $\gamma\coloneq\mathsf{min}\{\alpha<\delta:{T}_{\alpha+1}={T}_\alpha\}$.
Then
$$
{Q}\coloneq{T}_{\gamma}\qquad\text{and}\qquad
{D}\coloneq\bigcup_{\alpha\leqslant\gamma}({F}_\alpha\cup{G}_\alpha)
$$
work.

Indeed, ${Q}\subseteq{T}$ by (1) and (5),
and using (G) and (I) we have
$$
{D}\:\subseteq\:
\bigcup_{\alpha\leqslant\gamma}\mathsf{closure}({T}_\alpha,\SS^\omega)
\:\subseteq\:
\mathsf{closure}({T},\SS^\omega),
$$
so (a) is satisfied.
Condition (2) implies that
$\varphi\ll{F}_\alpha\cup{G}_\alpha\rr
\supseteq{T}_{\alpha}\setminus{T}_{\alpha+1}$
for all $\alpha<\delta$,
so using (5), (3), (1), and the choice of $\gamma$ we can write
$$
\varphi\ll{D}\rr\:=\:
\bigcup_{\alpha\leqslant\gamma}\varphi\ll{F}_\alpha\cup{G}_\alpha\rr
\:\supseteq\:
\bigcup_{\alpha\leqslant\gamma}({T}_{\alpha}\setminus{T}_{\alpha+1})
\:\supseteq\:
{T}_{0}\setminus{T}_{\gamma+1}\:=\:
{T}\setminus{T}_{\gamma}\:=\:{T}\setminus{Q},
$$
so (b) is satisfied.
Condition (2) implies $\varphi\ll{F}_\gamma\rr\cap{T}_{\gamma+1}=\varnothing$,
hence $\varphi\ll{F}_\gamma\rr\cap{T}_{\gamma}=\varnothing$, so it follows from (H) that
${T}_\gamma\cap\hedge({T}_\gamma)=\varnothing$.
Then ${T}_\gamma\subseteq\hcrown({T}_\gamma)$ by Remark~\ref{rem.psi}(f) and (6), so
(c) holds.
Similarly, we have $\varphi\ll{G}_\gamma\rr\cap{T}_{\gamma}=\varnothing$,
therefore $\varphi\ll{G}_\gamma\rr\cap\mathsf{closure}({T}_\gamma,\SS^\omega)=\varnothing$
because $\varphi\ll{G}_\gamma\rr$ is open in $\SS^\omega$.
Then (J) imply
$\mathsf{closure}({T}_\gamma,\SS^\omega)\cap\hbottom({T}_\gamma)=\varnothing$,
which means that (d) also holds.
It remains to prove that ${D}$ is a closed discrete set in $\SS^\omega$.

Suppose on the contrary that ${D}$ has an accumulation point ${r}$ in $\SS^\omega$;
so there is a sequence $\ll{p}_{i}\rr_{{i}\in\omega}$ in ${D}\setminus\{{r}\}$ that converges to ${r}$.
Then there is a sequence $\ll\alpha_{i}\rr_{{i}\in\omega}$ of ordinals
such that $\alpha_{i}\leqslant\gamma$ and ${p}_{i}\in{F}_{\alpha_i}\cup{G}_{\alpha_{i}}$ for all ${i}\in\omega$.
Passing to a subsequence, we may assume that one of two cases holds:

\emph{Case 1.} $\alpha_{i}=\alpha_{j}\ $ for all ${i},{j}\in\omega$.
This immediately leads to a contradiction with (4).

\emph{Case 2.} $\alpha_{i}<\alpha_{{i}+1}\ $ for all ${i}\in\omega$.
Passing to another subsequence, we may consider only two cases again:

\emph{Case 2.1.} ${p}_{i}\in{G}_{\alpha_{i}}\ $ for all ${i}\in\omega$.
It follows from (I) and Remark~\ref{rem.hat}(a) that

  \begin{itemize}
  \item [8.\,]
    $\dot{p}_{0}\in\bottom({T}_{\alpha_{0}})\nos\,$.
  \end{itemize}
Using Remark~\ref{lemm.C-lim} with (7) and passing to a subsequence once more, we may assume that
$\dot{p}_{i}\xrightarrow[{i}\to\infty]{{B}\downarrow}{\dot{r}}$
for some ${B}\subseteq{A}$.

\emph{Case 2.1.1.} ${B}={A}$.
Then $\dot{p}_1\sqsubseteq_{A}\dot{p}_0$, that is, $\dot{p}_1=\dot{p}_0$,
so $\dot{p}_{1}\in\bottom({T}_{\alpha_{0}})$ by (8),
whence ${p}_{1}\in\hbottom({T}_{\alpha_{0}})$ by Remark~\ref{rem.hat}(a).
Therefore, using (J), (5), (Case 2), (I), and (Case 2.1), we can write
$$
\varphi\ll{G}_{\alpha_0}\cup{F}_{\alpha_0}\rr
\:\supseteq\:
\varphi\ll{G}_{\alpha_0}\rr
\:\supseteq\:
\mathsf{closure}({T}_{\alpha_0},\SS^\omega)
\cap\hbottom({T}_{\alpha_0})
\:\supseteq\:
$$
$$
\mathsf{closure}({T}_{\alpha_1},\SS^\omega)
\cap\hbottom({T}_{\alpha_0})
\:\supseteq\:
{G}_{\alpha_1}\cap\hbottom({T}_{\alpha_0})
\:\ni\:{{p}_1},
$$
that is, $\varphi\ll{G}_{\alpha_0}\cup{F}_{\alpha_0}\rr$
is an open neighbourhood of ${p}_1$ in $\SS^\omega$.
Then, using (Case 2.1), (I), (5), (Case 2), and (2), we get
$$
{p}_1\:\in\:{G}_{\alpha_1}
\:\subseteq\:
\mathsf{closure}({T}_{\alpha_1},\SS^\omega)
\:\subseteq\:
\mathsf{closure}({T}_{\alpha_0+1},\SS^\omega)
\:=\:
$$
$$
\mathsf{closure}({T}_{\alpha_0}\setminus\varphi\ll{F}_{\alpha_0}\cup{G}_{\alpha_0}\rr,\SS^\omega)
\:\subseteq\:
\SS^\omega\setminus\varphi\ll{F}_{\alpha_0}\cup{G}_{\alpha_0}\rr
\:\subseteq\:
\SS^\omega\setminus\{{p}_1\},
$$
a contradiction.

\emph{Case 2.1.2.}
${B}\neq{A}$.
Since $\dot{p}_1\sqsubseteq_{B}\dot{p}_0$, we have
$\dot{p}_0\neq\dot{p}_1\sqsubseteq\dot{p}_0$
by Remark~\ref{rem.sqsubset}(a,c).
It follows from Remark~\ref{rem.psi}(i), (6), (7), and (8) that
${p}_{1}\notin\mathsf{closure}({T}_{\alpha_{0}},\SS^\omega)$.
This contradicts
$$
{p}_1\in{G}_{\alpha_1}\subseteq\mathsf{closure}({T}_{\alpha_1},\SS^\omega)
\subseteq\mathsf{closure}({T}_{\alpha_0},\SS^\omega),
$$
which follows from (Case 2.1), (I), (5), and (Case 2).

\emph{Case 2.2.}
${p}_{i}\in{F}_{\alpha_{i}}\ $ for all ${i}\in\omega$.
By (G), ${p}_{i}\in\mathsf{closure}    \big({T}_{\alpha_{i}}\cap\hedge({T}_{\alpha_{i}}),\SS^\omega\big)$ for all ${i}\in\omega$.
Using a countable local base at the point ${r}$ in $\SS^\omega$ and passing to a subsequence in $\ll{p}_{i}\rr_{{i}\in\omega}$, we may assume that there is a sequence $\ll{q}_{i}\rr_{{i}\in\omega}$ in $\SS^\omega$ that also converges to ${r}$ in $\SS^\omega$ and such that
${q}_{i}\in{T}_{\alpha_{i}}\cap\hedge({T}_{\alpha_{i}})$ for all ${i}\in\omega$.
Using Remark~\ref{lemm.C-lim} with (6) and passing to a subsequence again, we may assume that
$\dot{q}_{i}\xrightarrow[{i}\to\infty]{{C}\downarrow}\dot{r}$
for some ${C}\subseteq{A}$.
It follows from the choice of ${q}_1$, (5), and (Case 2) that
${q}_{1}\in{T}_{\alpha_{1}}\subseteq{T}_{\alpha_{0}}$;
also Remark~\ref{rem.sqsubset}(a), the choice of ${C}$, the choice of ${q}_0$, and Remark~\ref{rem.hat}(a) imply
$\dot{q}_{1}\sqsubseteq\dot{q}_{0}\in\edge({T}_{\alpha_{0}})$.
Then it follows from Remark~\ref{rem.psi}(g) that $\dot{q}_{1}\in\edge({T}_{\alpha_{0}})$,
so using ${q}_{1}\in{T}_{\alpha_{0}}$, (H), and Remark~\ref{rem.hat}(a), we get
$$
\varphi\ll{F}_{\alpha_{0}}\cup{G}_{\alpha_{0}}\rr
\:\supseteq\:
\varphi\ll{F}_{\alpha_{0}}\rr
\:\supseteq\:
{T}_{\alpha_{0}}\cap\hedge({T}_{\alpha_{0}})
\:\ni\:{q}_{1};
$$
that is, $\varphi\ll{F}_{\alpha_{0}}\cup{G}_{\alpha_{0}}\rr\ni{q}_{1}$.
Then it follows from the choice of ${q}_{1}$, (5), (Case 2), and (2) that
$$
{q}_{1}\:\in\:{T}_{\alpha_{1}}\:\subseteq\:{T}_{\alpha_{0}+1}
\:=\:{T}_{\alpha_{0}}\setminus\varphi\ll{F}_{\alpha_{0}}\cup{G}_{\alpha_{0}}\rr
\:\subseteq\:
\SS^\omega\setminus\varphi\ll{F}_{\alpha_{0}}\cup{G}_{\alpha_{0}}\rr
\:\subseteq\:
\SS^\omega\setminus\{{q}_{1}\},
$$
a contradiction.
\end{proof}

\section{Proofs of lemmas~\ref{lem.G} and~\ref{lem.F}}
\label{proof.of.lem.F}

\begin{proof}[\textbf{Proof of Lemma~\ref{lem.G}}]
Suppose that $\varphi$ is an NA for $\SS^\omega$, ${n}\in\omega$, and ${P}$ is a closed $\nos{n}$-broom on $\SS^\omega$.
Assume that $\mathsf{RIH}(\SS^\omega,\varphi,{n},{P})$ holds.
Let ${m}\in\omega$, ${m}>{n}$, $\varnothing\neq{A}\subseteq{m}$, ${\ZZ}$ be an ${A},{m}\nos$-cube, and ${T}\subseteq{P}\cap\breve\ZZ$ be an $({m}\!\setminus\!{A}\nos)$-broom on $\SS^\omega$.

\noindent
We must find a closed discrete set ${G}$ in $\SS^\omega$ such that

  \begin{itemize}
  \item [\textup{i.}\,]
    ${G}\subseteq
    \mathsf{closure}({T},\SS^\omega)\cap\hbottom({T})\quad\nos\,$ and
  \item [\textup{j.}\,]
    $\varphi\ll{G}\rr\supseteq
    \mathsf{closure}({T},\SS^\omega)\cap\hbottom({T})\nos\,$.
  \end{itemize}

Put ${W}_{x}\coloneq\{{x}\}$ for all ${x}\in{\ZZ}$.
Then $\bottom({T})=\bigcup\{{W}_{x}:{x}\in\bottom({T})\}$,
so $\hbottom({T})=\bigcup\{\breve{W}_{x}:{x}\in\bottom({T})\}$
by Remark~\ref{rem.hat}(a).
We have
$$
\mathsf{closure}({T},\SS^\omega)\cap\hbottom({T})
\:=\:
\bigcup_{{x}\in\bottom({T})}\big(\mathsf{closure}({T},\SS^\omega)\cap\breve{W}_{x}\big)
\:=\:\bigcup_{{x}\in\bottom({T})}{Q}_{x},
$$
where ${Q}_{x}\coloneq\mathsf{closure}({T},\SS^\omega)\cap\breve{W}_{x}$
for all ${x}\in\bottom({T})$.
Each $\breve{W}_{x}$ is an $\nos{A}$-broom on $\SS^\omega$,
so each ${Q}_{x}$ is also an $\nos{A}$-broom by Remark~\ref{rem.brooms}(c).
Since ${T}$ is an $({m}\!\setminus\!{A}\nos)$-broom,
it follows from Remark~\ref{rem.brooms}(a,c) that
each ${Q}_{x}$ is also an $({m}\!\setminus\!{A}\nos)$-broom.
Therefore each ${Q}_{x}$ is an ${m}\nos$-broom
by Remark~\ref{rem.brooms}(d).
Also each ${Q}_{x}$ is a closed subset of $\SS^\omega$
and ${Q}_{x}\subseteq{P}$ because ${T}\subseteq{P}$ and ${P}$ is closed in $\SS^\omega$.
Then by $\mathsf{RIH}(\SS^\omega,\varphi,{n},{P})$,
for each ${x}\in\bottom({T})$,
there exists a closed discrete set ${G}_{x}$ in the subspace ${Q}_{x}$ such that $\varphi\ll{G}_{x}\rr\supseteq{Q}_{x}$.
Then each ${G}_{x}$  is a closed discrete set in $\SS^\omega$.

We prove that the set ${G}\coloneq\bigcup_{{x}\in\bottom({T})}{G}_{x}$ satisfies the required conditions.
We have
$$
{G}
\:=\:
\bigcup_{{x}\in\bottom({T})}{G}_{x}
\:\subseteq\:
\bigcup_{{x}\in\bottom({T})}{Q}_{x}
\:=\:
\mathsf{closure}({T},\SS^\omega)\cap\hbottom({T}),
$$
so (i) is satisfied.
We also have
$$
\varphi\ll{G}\rr
\:=\:
\bigcup_{{x}\in\bottom({T})}\varphi\ll{G}_{x}\rr
\:\supseteq\:
\bigcup_{{x}\in\bottom({T})}{Q}_{x}
\:=\:
\mathsf{closure}({T},\SS^\omega)\cap\hbottom({T}),
$$
so (j) is also satisfied.
It remains to show that ${G}$ is a closed discrete set in $\SS^\omega$.

Suppose on the contrary that ${G}$ has an accumulation point ${r}$ in $\SS^\omega$;
so there is a sequence $\ll{p}_{i}\rr_{{i}\in\omega}$ in ${G}\setminus\{{r}\}$ that converges to ${r}$.
Then there is a sequence $\ll{x}_{i}\rr_{{i}\in\omega}$ in $\bottom({T})$
such that ${p}_{i}\in{G}_{{x}_{i}}$ for all ${i}\in\omega$.
For all ${i}\in\omega,$ ${p}_{i}\in{G}\subseteq\hbottom({T})\subseteq\breve\ZZ$
by (i) and Remark~\ref{rem.hat}(b),
so using Remark~\ref{lemm.C-lim} and passing to a subsequence, we may assume that
$\dot{p}_{i}\xrightarrow[{i}\to\infty]{{B}\downarrow}{\dot{r}}$
for some ${B}\subseteq{A}$.
For all ${i}\in\omega\setminus\{0\}$,
we have $\dot{p}_{i}\sqsubseteq\dot{p}_0$
by Remark~\ref{rem.sqsubset}(a,b),
so, for all ${i}\in\omega\setminus\{0\}$, $\dot{p}_{i}=\dot{p}_0$
because $\dot{p}_{i},\dot{p}_0\in\bottom({T})$ by Remark~\ref{rem.hat}(a) and (i).
We have
$$
{p}_{i}
\:\in\:
{G}_{{x}_{i}}
\:\subseteq\:
{Q}_{{x}_{i}}
\:\subseteq\:
\breve{W}_{{x}_{i}}
$$
for all ${i}\in\omega$,
therefore ${x}_{i}=\dot{p}_{i}=\dot{p}_0={x}_{0}$ for all ${i}\in\omega$.
So ${p}_{i}\in{G}_{{x}_{i}}={G}_{{x}_{0}}$ for all ${i}\in\omega$,
hence ${G}_{{x}_{0}}$ has an accumulation point in $\SS^\omega$,
a contradiction.
\end{proof}
\medskip

To prove Lemma~\ref{lem.F}, we need additional notation and two lemmas, which we prove in the next section.

\begin{nota}\label{nota.harpoon}
Suppose that ${P}\subseteq\prod_{{i}\in{A}}{X}_{i}$ and ${B}\subseteq{A}$.
Then:

  \begin{itemize}
  \item [\ding{46}\:]
    ${P}\,|\,{B}\:\coloneq\:
    \{{p}\res{B}:{p}\in{P}\}\nos$.
  \end{itemize}
\end{nota}

\begin{nota}
Suppose that ${m}\in\omega$, $\varnothing\neq{A}\subseteq{m}$, ${\ZZ}$ is an ${A},{m}\nos$-cube, and ${T}\subseteq\breve\ZZ$. Let ${B}\subseteq{A}$, ${w}\in{\ZZ}\,|\,{B}$, and ${x}\in{\ZZ}$.
Then:

  \begin{itemize}
  \item [\ding{46}\:]
    $\plane({B},{w})\:\coloneq\:
    \{{y}\in{\ZZ}:{y}\res{B}={w}\}\nos\,$.
  \item [\ding{46}\:]
    ${x}\nos\,$ is ${B}\text{\emph{-mobile}}\iff
    {x}\in\edge({T})\quad\text{and}\quad
    \exists{y}\in\edge({T})\:[{y}\sqsubseteq_{B}{x}]\nos\,$.
  \item [\ding{46}\:]
    $\facet({B},{w})\:\coloneq\:
    \{{x}\in\edge({T})\cap\plane({B},{w}):
    {x}\text{ is }{B}\text{-mobile and }
    \neg\exists{E}\subsetneqq{B}\:[{x}\text{ is }{E}\text{-mobile}]\}\nos\,$.
  \end{itemize}
\end{nota}
\noindent
Note that the above notations depend on ${m}$, ${A}$, ${\ZZ}$, and ${T}$, but we write, for example, $\facet({B},{w})$ instead of $\facet_{{A},{m},{\ZZ},{T}}({B},{w})$ for better readability.

\begin{rema}\label{rem.plane}
Suppose that ${m}\in\omega$, $\varnothing\neq{A}\subseteq{m}$, ${\ZZ}$ is an ${A},{m}\nos$-cube, and ${T}\subseteq\breve\ZZ$. Let ${B}\subseteq{A}$ and ${w}\in{\ZZ}\,|\,{B}$.
Then:

  \begin{itemize}
  \item [\textup{a.}\,]
    $\hplane({B},{w})\nos\,$ is a closed subset of $\SS^\omega$.
  \item [\textup{b.}\,]
    $\hplane({B},{w})\nos\,$ is a ${B}\nos$-broom on $\SS^\omega$.\hfill$\qed$%
  \end{itemize}%
\end{rema}%
\smallskip

\begin{lemm}\label{lem.facet}
Suppose that ${m}\in\omega$, $\varnothing\neq{A}\subseteq{m}$, ${\ZZ}$ is an ${A},{m}\nos$-cube, and ${T}\subseteq\breve\ZZ$.
Then

  \begin{itemize}
  \item [\ding{43}\,]
    $\hedge({T})\:\subseteq\:\bigcup\{\hfacet({B},{w}):
    \varnothing\neq{B}\subseteq{A}\text{ and } {w}\in{\ZZ}\,|\,{B}\}\nos\,$.
  \end{itemize}%
\end{lemm}%
\smallskip

\begin{lemm}\label{lem.F(C,w)}
Suppose that $\varphi$ is an NA for $\SS^\omega$, ${n}\in\omega$, ${P}$ is a closed $\nos{n}$-broom on $\SS^\omega$, and ${l}\in\omega\setminus\{0\}$.
Assume that $\mathsf{IH}(\varphi,{n},{P},{l})$ holds.
Let ${m}\in\omega$, ${A}\in[{m}]^l$, ${\ZZ}$ be an ${A},{m}\nos$-cube,
$\varnothing\neq{B}\subseteq{A}$, ${w}\in{\ZZ}\,|\,{B}$, and ${T}\subseteq{P}\cap\breve\ZZ$ be an $({m}\!\setminus\!{A}\nos)$-broom on $\SS^\omega$  such that $\mathsf{rank}({t},\varphi)={m}$ for all ${t}\in{T}$.

Then there is a closed discrete set ${H}$ in $\SS^\omega$ such that

  \begin{itemize}
  \item [\textup{\textbf{g}}.\,]
    ${H}\subseteq\mathsf{closure}
    \big({T}\cap\hfacet({B},{w}),\SS^\omega\big)\quad\nos\,$ and
  \item [\textup{\textbf{h}}.\,]
    $\varphi\ll{H}\rr\supseteq{T}\cap\hfacet({B},{w})\nos\,$.
  \end{itemize}
\end{lemm}%
\smallskip

\begin{proof}[\textbf{Proof of Lemma~\ref{lem.F}}]
Suppose that $\varphi$ is an NA for $\SS^\omega$, ${n}\in\omega$, ${P}$ is a closed $\nos{n}$-broom on $\SS^\omega$, and ${l}\in\omega\setminus\{0\}$.
Assume that $\mathsf{IH}(\varphi,{n},{P},{l})$ holds.
Let ${m}\in\omega$, ${A}\in[{m}]^l$, ${\ZZ}$ be an ${A},{m}\nos$-cube, and ${T}\subseteq{P}\cap\breve\ZZ$ be an $({m}\!\setminus\!{A}\nos)$-broom on $\SS^\omega$
such that $\mathsf{rank}({t},\varphi)={m}$ for all ${t}\in{T}$.

\noindent
We must find a closed discrete set ${F}$ in $\SS^\omega$ such that

  \begin{itemize}
  \item [\textup{g.}\,]
    ${F}\subseteq\mathsf{closure}
    \big({T}\cap\hedge({T}),\SS^\omega\big)\quad\nos\,$ and
  \item [\textup{h.}\,]
    $\varphi\ll{F}\rr\supseteq{T}\cap\hedge({T})\nos\,$.
  \end{itemize}

Using Lemma~\ref{lem.F(C,w)}, for each nonempty ${B}\subseteq{A}$ and ${w}\in{\ZZ}\,|\,{B}$, let ${H}({B},{w})$ be a closed discrete set in $\SS^\omega$ such that

  \begin{itemize}
  \item [\textup{\textbf{G}}.\,]
    ${H}({B},{w})\subseteq\mathsf{closure}
    \big({T}\cap\hfacet({B},{w}),\SS^\omega\big)\quad\nos\,$ and
  \item [\textup{\textbf{H}}.\,]
    $\varphi\ll{H}({B},{w})\rr\supseteq{T}\cap\hfacet({B},{w})\nos\,$.
  \end{itemize}
Then the set
$$\textstyle
{F}\:\coloneq\:\bigcup\{{H}({B},{w}):\varnothing\neq{B}\subseteq{A}\text{ and } {w}\in{\ZZ}\,|\,{B}\}
$$
satisfies the required conditions.
We have $\hfacet({B},{w})\subseteq\hedge({T})$ by Remark~\ref{rem.hat}(b),
so each ${H}({B},{w})$
is a subset of $\mathsf{closure}\big({T}\cap\hedge({T}),\SS^\omega\big)$,
from which (g) follows.
Using (\textup{\textbf{H}}) and Lemma~\ref{lem.facet}, we can write
$$\textstyle
\varphi\ll{F}\rr
\:=\:
\bigcup\{\varphi\ll{H}({B},{w})\rr:\varnothing\neq{B}\subseteq{A}\text{ and } {w}\in{\ZZ}\,|\,{B}\}
\:\supseteq\:
$$
$$\textstyle
\bigcup\{{T}\cap\hfacet({B},{w}):\varnothing\neq{B}\subseteq{A}\text{ and } {w}\in{\ZZ}\,|\,{B}\}
\:\supseteq\:
$$
$$\textstyle
{T}\cap\bigcup\{\hfacet({B},{w}):\varnothing\neq{B}\subseteq{A}\text{ and } {w}\in{\ZZ}\,|\,{B}\}
\:\supseteq\:
{T}\cap\hedge({T}),
$$
so (h) is also satisfied.
It remains to prove that ${F}$ is a closed discrete set in $\SS^\omega$.

Suppose on the contrary that ${F}$ has an accumulation point ${r}$ in $\SS^\omega$;
so there is a sequence $\ll{p}_{i}\rr_{{i}\in\omega}$ in ${F}\setminus\{{r}\}$ that converges to ${r}$.
Then there are sequences $\ll{B}_{i}\rr_{{i}\in\omega}$
and $\ll{w}_{i}\rr_{{i}\in\omega}$ such that $\varnothing\neq{B}_{i}\subseteq{A}$,
${w}_{i}\in{\ZZ}\,|\,{B}_{i}$,
and
${p}_{i}\in{H}({B}_{i},{w}_{i})$ for all ${i}\in\omega$.
Because ${A}$ is finite, passing to a subsequence, we may assume that
there is nonempty ${B}\subseteq{A}$ such that
${B}_{i}={B}$ for all ${i}\in\omega$.
It follows form (\textup{\textbf{G}}) and Remark~\ref{rem.hat}(c) that ${F}\subseteq\breve\ZZ$,
so using Remark~\ref{lemm.C-lim} and passing to another subsequence, we may assume that
$\dot{p}_{i}\xrightarrow[{i}\to\infty]{{C}\downarrow}{\dot{r}}$
for some ${C}\subseteq{A}$.

Using (\textup{\textbf{G}}), for each ${i}\in\omega$,
we can find a sequence $\ll{q}_{{i},{j}}\rr_{{j}\in\omega}$
in ${T}\cap\hfacet({B},{w}_{i})$ that converges to ${p}_{i}$ in $\SS^\omega$.
Then, for each ${i}\in\omega$, using Remark~\ref{lemm.C-lim}
and passing to a subsequence in the sequence $\ll{q}_{{i},{j}}\rr_{{j}\in\omega}$,
we may assume that $\dot{q}_{{i},{j}}\xrightarrow[{j}\to\infty]{{D}_{i}\downarrow}{\dot{p}_{i}}$
for some ${D}_{i}\subseteq{A}$.
Passing to another subsequence, now again in the sequence $\ll{p}_{i}\rr_{{i}\in\omega}$,  we may assume that there is ${D}\subseteq{A}$ such that
${D}_{i}={D}$ for all ${i}\in\omega$.

We have $\dot{p}_{1}\sqsubseteq_{C}\dot{p}_{0}\sqsubseteq_{D}\dot{q}_{0,0}$
by the choice of ${C}$, ${D}$, and ${D}_{0}$,
so $\dot{p}_{1}\sqsubseteq_{{C}\cap{D}}\dot{q}_{0,0}$ by Remark~\ref{rem.sqsubset}(b).
Then, since
$\dot{q}_{1,{j}}\xrightarrow[{j}\to\infty]{{D}\downarrow}{\dot{p}_1}$,
it follows from Remark~\ref{rem.1.25} that
there is ${{j}'}\in\omega$ such that
$\dot{q}_{1,{{j}'}}\sqsubseteq_{{C}\cap{D}}\dot{q}_{0,0}$.
Since ${q}_{1,{{j}'}}\in\hfacet({B},{w}_{1})$,
we have $\dot{q}_{1,{{j}'}}\in\facet({B},{w}_{1})$ by Remark~\ref{rem.hat}(a),
therefore $\dot{q}_{1,{{j}'}}$ is ${B}\nos$-mobile,
so there is ${y}\in\edge({T})$
such that ${y}\sqsubseteq_{B}\dot{q}_{1,{{j}'}}$.
Then $\dot{q}_{1,{{j}'}}\sqsubseteq_{{C}\cap{D}}\dot{q}_{0,0}$
implies ${y}\sqsubseteq_{{B}\cap{C}\cap{D}}\dot{q}_{0,0}$
by Remark~\ref{rem.sqsubset}(b).
Therefore $\dot{q}_{0,0}$ is $({B}\cap{C}\cap{D})\nos$-mobile
because ${y}\in\edge({T})$
and $\dot{q}_{0,0}\in\facet({B},{w}_{0})\subseteq\edge({T})$
(since ${q}_{0,0}\in\hfacet({B},{w}_{0})$).
Now, since $\dot{q}_{0,0}\in\facet({B},{w}_{0})$, it follows that
there is no ${E}\subsetneqq{B}$ such that $\dot{q}_{0,0}$ is ${E}\nos$-mobile.
But $\dot{q}_{0,0}$ is $({B}\cap{C}\cap{D})\nos$-mobile and
$({B}\cap{C}\cap{D})\subseteq{B}$,
therefore $({B}\cap{C}\cap{D})={B}$,
which implies ${B}\subseteq{C}$.

We have $\dot{p}_{i}\xrightarrow[{i}\to\infty]{{C}\downarrow}{\dot{r}}$
by the choice of ${C}$,
therefore $\dot{p}_{i}\res{C}=\dot{r}\res{C}$ for all ${i}\in\omega$, so ${B}\subseteq{C}$ implies

  \begin{itemize}
  \item [1.\,]
    $\dot{p}_{i}\res{B}=\dot{r}\res{B}\quad\nos\,$
    for all ${i}\in\omega$.
  \end{itemize}
Also, by the choice of ${B}$, ${B}_{i}$, and ${w}_{i}$,
by (\textup{\textbf{G}}), by the definition of $\facet({B},{w}_{i})$ with Remark~\ref{rem.hat}(b),
and by Remark~\ref{rem.plane}(a),
for all ${i}\in\omega$,
$$
{p}_{i}
\:\in\:
{H}({B},{w}_{i})
\:\subseteq\:
\mathsf{closure}\big(\hfacet({B},{w}_{i}),\SS^\omega\big)
\:\subseteq\:
\mathsf{closure}\big(\hplane({B},{w}_{i}),\SS^\omega\big)
\:=\:
\hplane({B},{w}_{i}).
$$
Then, for all ${i}\in\omega$, $\dot{p}_{i}\in\plane({B},{w}_{i})$ by Remark~\ref{rem.hat}(a),
so $\dot{p}_{i}\res{B}={w}_{i}$,
and then ${w}_{i}=\dot{r}\res{B}$ by (1).
Therefore, for all ${i}\in\omega$, we have ${w}_{i}={w}_{0}$,
hence ${p}_{i}\in{H}({B},{w}_{i})={H}({B}_{0},{w}_{0})$.
It follows that ${r}$ is an accumulation point of the set ${H}({B}_{0},{w}_{0})$,
a contradiction.
\end{proof}

\section{Proofs of Lemmas~\ref{lem.facet} and~\ref{lem.F(C,w)}}

\begin{proof}[\textbf{Proof of Lemma~\ref{lem.facet}}]
Suppose that ${m}\in\omega$, $\varnothing\neq{A}\subseteq{m}$, ${\ZZ}$ is an ${A},{m}\nos$-cube, and ${T}\subseteq\breve\ZZ$. We must show that

  \begin{itemize}
  \item [\ding{43}\,]
    $\hedge({T})\:\subseteq\:\bigcup\{\hfacet({B},{w}):
    \varnothing\neq{B}\subseteq{A}\text{ and } {w}\in{\ZZ}\,|\,{B}\}\nos\,$.
  \end{itemize}

Assume that ${p}\in\hedge({T})$; then
$\dot{p}\in\edge({T})$ by Remark~\ref{rem.hat}(a).
For each ${y}\in{\ZZ}$ with ${y}\sqsubseteq{\dot{p}}$,
there is the unique ${B}({y})\subseteq{A}$
such that ${y}\sqsubseteq_{{B}({y})}{\dot{p}}$.
Let ${Y}\coloneq\{{y}\in\edge({T}):{y}\sqsubseteq{\dot{p}}\}$.
This set is not empty because ${\dot{p}}\in{Y}$,
so there is ${y}'\in{Y}$ such that
$|{B}({y}')|\leqslant|{B}({y})|$ for all ${y}\in{Y}$.
Let ${w}\coloneq{\dot{p}}\res{B}({y}')\in{\ZZ}\,|\,{B}({y}')$.
We have ${\dot{p}}\in\plane({B}({y}'),{w})$
and ${\dot{p}}$ is ${B}({y}')\nos$-mobile
because ${y}'\in\edge({T})$ and ${y}'\sqsubseteq_{{B}({y}')}{\dot{p}}$.
Then, by the choice of ${y}'$,
there is no ${E}\subsetneqq{B}({y}')$
such that ${\dot{p}}$ is ${E}\nos$-mobile,
therefore ${\dot{p}}\in\facet\big({B}({y}'),{w}\big)$,
so ${p}\in\hfacet\big({B}({y}'),{w}\big)$.
It remains to show that ${B}({y}')\neq\varnothing$.

Suppose on the contrary that ${B}({y}')$ is empty,
so that ${y}'\sqsubseteq_\varnothing{\dot{p}}$.
Since ${y}'\in\edge({T})$, we have
${y}'\in\mathsf{closure}\big(\crown({T}),{\ZZ}\big)$,
so by Remark~\ref{rem.sqsubset}(d) there is ${z}\in\crown({T})$
such that ${z}\sqsubseteq_\varnothing{\dot{p}}$.
Then ${\dot{p}}\in\crown({T})$ by Remark~\ref{rem.psi}(b)
and Remark~\ref{rem.hat}(a),
hence ${\dot{p}}\notin\edge({T})$ by Remark~\ref{rem.psi}(c), a contradiction.
\end{proof}
\medskip

\begin{proof}[\textbf{Proof of Lemma~\ref{lem.F(C,w)}}]
Suppose that $\varphi$ is an NA for $\SS^\omega$, ${n}\in\omega$, ${P}$ is a closed $\nos{n}$-broom on $\SS^\omega$, and ${l}\in\omega\setminus\{0\}$.
Assume that $\mathsf{IH}(\varphi,{n},{P},{l})$ holds.
Let ${m}\in\omega$, ${A}\in[{m}]^l$, ${\ZZ}$ be an ${A},{m}\nos$-cube,
$\varnothing\neq{B}\subseteq{A}$, ${w}\in{\ZZ}\,|\,{B}$, and ${T}\subseteq{P}\cap\breve\ZZ$ be an $({m}\!\setminus\!{A}\nos)$-broom on $\SS^\omega$  such that $\mathsf{rank}({t},\varphi)={m}$ for all ${t}\in{T}$.

\noindent
We must find a closed discrete set ${H}$ in $\SS^\omega$ such that

  \begin{itemize}
  \item [\textup{\textbf{g}}.\,]
    ${H}\subseteq\mathsf{closure}
    \big({T}\cap\hfacet({B},{w}),\SS^\omega\big)\quad\nos\,$ and
  \item [\textup{\textbf{h}}.\,]
    $\varphi\ll{H}\rr\supseteq{T}\cap\hfacet({B},{w})\nos\,$.
  \end{itemize}

Put ${R}\coloneq\mathsf{closure}\big({T}\cap\hfacet({B},{w}),\SS^\omega\big)$.
We have ${R}\subseteq{P}$ because ${T}\subseteq{P}$ and ${P}$ is closed in $\SS^\omega$.
By Remark~\ref{rem.hat}(b) and Remark~\ref{rem.plane}(a),
$$
{R}
\:\subseteq\:
\mathsf{closure}\big(\hfacet({B},{w}),\SS^\omega\big)
\:\subseteq\:
\mathsf{closure}\big(\hplane({B},{w}),\SS^\omega\big)
\:=\:
\hplane({B},{w}),
$$
so ${R}$ is a ${B}\nos$-broom on $\SS^\omega$ by Remark~\ref{rem.brooms}(c) and  Remark~\ref{rem.plane}(b).
Also ${R}\subseteq\mathsf{closure}({T},\SS^\omega)$,
so ${R}$ is an $({m}\!\setminus\!{A}\nos)$-broom by Remark~\ref{rem.brooms}(a,c).
Then ${R}$ is a $\big({B}\cup({m}\!\setminus\!{A})\big)\nos$-broom by Remark~\ref{rem.brooms}(d).
Put ${C}\coloneq{A}\setminus{B}$;
then ${m}\setminus{C}={m}\setminus({A}\!\setminus\!{B})
\subseteq{B}\cup({m}\!\setminus\!{A})$,
so ${R}$ is an $({m}\!\setminus\!{C}\nos)$-broom
by Remark~\ref{rem.brooms}(b). We have $|{C}|<|{A}|={l}$
because $\varnothing\neq{B}\subseteq{A}$,
so, by $\mathsf{IH}(\varphi,{n},{P},{l})$,
$\text{\ding{74}}(\varphi,{n},{P},|{C}|)$ holds.

Now, ${m}\in\omega$, ${C}\in[{m}]^{|{C}|}$, and ${R}\subseteq{P}$ is a closed $({m}\!\setminus\!{C}\nos)$-broom on $\SS^\omega$,
so by $\text{\ding{74}}(\varphi,{n},{P},|{C}|)$
there exists a closed discrete set ${H}$ in $\SS^\omega$ such that

  \begin{itemize}
  \item [\ding{40}\,]
    ${H}\subseteq{R}\quad\nos\,$ and
  \item [\ding{40}\,]
    $\varphi\ll{H}\rr\supseteq
    \{{r}\in{R}:\mathsf{rank}({r},\varphi)={m}\}\nos\,$.
  \end{itemize}
Then (\textup{\textbf{g}}) is satisfied because ${H}\subseteq{R}$.
At last, since $\mathsf{rank}({r},\varphi)={m}$ for all ${r}\in{T}$,
we have
$$
\nos\varphi\ll{H}\rr
\:\supseteq\:
\big\{{r}\in\mathsf{closure}\big({T}\cap\hfacet({B},{w}),\SS^\omega\big)
:\mathsf{rank}({r},\varphi)={m}\big\}
\:\supseteq\:
$$
$$
\big\{{r}\in{T}\cap\hfacet({B},{w}):\mathsf{rank}({r},\varphi)={m}\big\}
\:\supseteq\:
{T}\cap\hfacet({B},{w}),
$$
so (\textup{\textbf{h}}) holds.
\end{proof}


\end{document}